\documentclass[12pt]{article}
\usepackage{amsmath}
\usepackage{amssymb}
\usepackage{amsthm}
\usepackage{amsfonts}
\usepackage{graphicx}

\usepackage{bbm,epsfig,graphics,epic,color,rotating,color}
\numberwithin{equation}{section}

\newcommand{\Z}{{\mathbb Z}}

\newcommand{\wt}[1]{\widetilde{#1}}

\def\bQ{\mathbf Q}\def\vphi{\varphi} 
\def\rb{{\rm b}} \def\rC{{\rm C}}
\def\rd{{\rm d}} \def\rD{{\rm D}} \def\rT{{\rm T}}
\def\rx{{\rm x}} 
\def\gam{\gamma}  
 \def\otau{\overline\tau}

\def\uy{{\underline y}} \def\bbZ{{\mathbb Z}} 
\def\eps{\epsilon}
\def\ups{\upsilon}
\def\om{\omega}

\def\Lam{{\Lambda}}
\def\lam{{\lambda}}
\def\ulam{\underline\lam}
\def\umu{{\underline\mu}}
\def\unn{{\underline n}}
\def\ue{{\underline e}}
\def\urho{{\underline \rho}}
\def\urhomu{\underline{\rho\iota}}
\def\diy{\displaystyle}
\def\beq{\begin{equation}}
\def\eeq{\end{equation}}

\theoremstyle{definition}

\begin{document}\def\cl{\centerline}
\def\Lam{{\Lambda}}
\def\lam{{\lambda}}
\def\ulam{\underline\lam}
\def\umu{{\underline\mu}}
\def\unn{{\underline n}} \def\uNn{{\underline{n_{}}}}
\def\ue{{\underline e}}
\def\urho{{\underline \rho}}
\def\urhomu{\underline{\rho\iota}}
\def\diy{\displaystyle}
\def\beq{\begin{equation}}
\def\eeq{\end{equation}}

\theoremstyle{definition}

\def\obeta{\overline\beta} \def\bbR{{\mathbb R}} \def\unu{{\underline\nu}}
\def\ka{{\kappa}} \def\ups{{\upsilon}}
\def\urho{{\underline\rho}} \def\usi{{\underline\sigma}}  \def\rp{{\rm p}} 
\def\bka{{\mbox{\boldmath$\kappa$}}} \def\bsi{{\mbox{\boldmath$\sigma$}}}
\def\bpi{{\mbox{\boldmath$\pi$}}} \def\bXi{{\mbox{\boldmath$\Xi$}}}
\def\Gam{{\Gamma}} \def \ux{{\underline x}}  \def\uy{{\underline y}}  
\def\ualp{{\underline\alpha}}

\def\BA{{\mathbf A}} \def\BL{{\mathbf L}}  \def\BI{{\mathbf I}}
\def\BM{{\mathbf M}}  \def\BR{{\mathbf R}}  \def\BW{{\mathbf W}}

\def\bbC{{\mathbb C}} \def\bbF{{\mathbb F}} \def\bbG{{\mathbb G}} \def\bbX{{\mathbb X}}

 \def\cF{{\mathcal F}}  \def\cZ{{\mathcal Z}}
 
\def\rC{{\rm C}} \def\rD{{\rm D}}  \def\rL{{\rm L}}
\def\rP{{\rm P}} \def\rI{{\rm I}}

\def\tA{{\tt A}} \def\tE{{\tt E}} \def\tF{{\tt F}} \def\tI{{\tt I}} \def\tL{{\tt L}} \def\tM{{\tt M}}
\def\tP{{\tt P}} \def\tQ{{\tt Q}} \def\tR{{\tt R}} \def\tS{{\tt S}}

\def\cmp{{\complement}}

\cl{\LARGE{\bf Models of Markov processes}}
\vskip .5 truecm
\cl{\LARGE{\bf with a random transition mechanism}}
\vskip 1 truecm

\cl{\Large{\bf Y. BelopolÕskaya$^1$,\ Y. Suhov$^{2-4}$}} 

\date{\today}

\vskip .5 truecm

{\bf Abstract.} The paper deals with a certain class of random evolutions. We develop a construction 
that yields 
an invariant measure for a continuous-time Markov process with random transitions. The approach 
is based on a 
particular way of constructing the combined process, where the generator 
is defined as a sum of two terms: one responsible for the evolution of the environment and the
second representing generators 
of processes with a given state of environment. (The two operators are not assumed to commute.) 
The presentation includes fragments of a general theory and pays a particular attention 
to several types of examples: (1) a queueing system with a random change of parameters (including
a Jackson network and, as a special case: a single-server queue with a diffusive
behavior of arrival and service rates), (2) a simple-exclusion model in presence of a special `heavy` particle,
(3) a diffusion with drift-switching, and (4) a diffusion with a randomly 
diffusion-type
varying diffusion coefficient (including a modification of the Heston random volatility model).   

\vskip .5 truecm



\section{Introduction}\label{sec:intro} 

This paper presents a construction of Markovian models of random processes in
a random environment. The problem can be stated in the following form. Suppose we are 
given a family of `basic`\\  

--------------------------------

$^1$St Petersburg State University of Architecture and Civil Engineering, 190005 RF

$^2$IITP, Moscow-127994, RF; $^5$DPMMS, University of Cambridge, CB3 0WB UK; 
$^3$Math Department, PSU, State College, PA16802 USA; 
$^4$IHES, Bures-sur-Yvette, 91440 France;

 E-mail: yms@statslab.cam.ac.uk; ims14@psu.edu 

\newpage

\noindent continuous-time Markov processes (MPs)  ${\wt X}^{(z)}(t)$ on a state space $\bbX$ 
where $z$ is a parameter
describing a `state of environment` (SE) which can be varied within a space $\cZ$. Assume that
$\forall$ $z\in\cZ$, process ${\wt X}^{(z)}$ has an invariant measure (IM) $\nu^{(z)}$. Further, 
suppose that we have
an MP ${\wt Z}(t)$ on $\cZ$ with an IM $\ups$. Is it possible to construct a `combined` 
(two-component) MP $(Z(t),X(t))$ on 
$\cZ\times\bbX$ in a `universal` manner, allowing a natural interpretation of an MP in a random 
environment? In this paper we put forward such a construction (under certain limitations); 
a feature of our constructions is a mutual 
impact of an SE and a state of a basic process upon each other. In other words, we describe a 
general mechanism where the state of a basic process influences a change of an environment
which in turn leads to a new transition rule for the process. Such an approach should be
compared with a body of work on random  walks in random environment (cf., e.g., \cite{Ze}) where
an environment is randomly chosen but kept fixed throughout the time dynamics.

An example of a system to which our construction can be applied is a Jackson network (JN) model;
see \cite{J1}, \cite{J2}. Here the environment can be identified with a triple $(\ulam ,\umu ,\tP)$ where
$\ulam$ and $\umu$ are vectors of arrival and service rates, and $\tP$ is a routing matrix.
The IM is a product of geometric distributions identified in terms of $\ulam$, $\umu$ and $\tP$. Making 
triple $(\ulam ,\umu ,\tP)$ dependent on parameter $z$ varying randomly within a finite or countable set $\cZ$
leads to a number of interesting applications, viz., networks with `distinguished  customers`;
cf. \cite{GPSY} and \cite{YS}. Our construction for this example is carried in Section \ref{sec:2};
it leads to a `modified` product-form IM for the combined MP and -- after the normalization -- to
a combined equilibrium probability distribution (EPD). 


In Section \ref{sec:3} we discuss another instructive example: a simple exclusion model on 
a cubic lattice $\bbZ^d$ interacting with a special `heavy` particle. Here again, the
poduct-form for an EPD of  
a combined MP is preserved in the course of the cnstruction.     

In Section \ref{sec:4}  we turn to general continuous-time MPs. First,  we cover a general class 
of jump MPs. Then a general form of our construction is provided, explaining the main 
mechanism behind it.  In Sections \ref{sec:5} and \ref{sec:6} we focus upon examples where 
one or both components
are diffusions. This includes queues where parameters exhibit 
a diffusion-type behavior as well as examples of a diffusion with a jump-like change of 
the drift coefficient. In particular, in Section \ref{sec:6} we discuss several models related
to mathematical finances: here basic 
MPs are Ornstein--Uhlenbeck diffusion processes with varying volatility.
 
A number of results in this paper are stated in the form of weak invariance equations (WIEs)
involving a measure and an operator    
acting on functions of a combined  state (in a simplified form -- a collection of combined jump rates).
In cases where we are able to assert the existence/uniqueness
of a combined MP, the WIE implies a genuine character of an IM; in these cases 
we refer to the invariance property directly. For convenience,
the WIE is stated anew in each specific (or general) context it is used.    
It has to be stressed again that our results depend on a special choice of the transition mechanism 
(inherited from Ref. \cite{GPSY})
where the SE and the basic MPs have a particular influence on each other. 

Throughout the paper, we repeatedly quote books \cite{EK}, \cite{KS}, \cite{Lig1} and \cite{Lig3} for various 
general results on
MPs and their generators. In fact, our models can be considered as examples of random 
evolution models considered in \cite{EK}, Chapter 12, although we focus here on different aspects
of behavior. Also, in contrast to  \cite{EK}, we -- as was said above -- include a mutual impact 
of the environment and 
basic MPs upon each other, albeit in a rather specific form.  (To make a comparison: the 
corresponding terms used in \cite{EK} are {\it driving} and {\it driven} processes.) Technically, in the 
provided assertions we attempt to stick to minimal assumptions (as a rule, mentioning them in passing).
But a mathematically minded reader should pay attention to remarks where we build links with 
more advanced conclusions by using general results from \cite{EK}, \cite{KS}, \cite{Lig1} and \cite{Lig3}.  


\section{A Jackson network in a random environment}\label{sec:2}

{\bf 2.A. Open Jackson networks.} \ In this section we use a Jackson job-shop 
network as a background model of a basic MP. The JN model is defined by the following ingredients
\cite{J1,J2}.

\vskip .5 truecm
\def\oR{{\overline R}}
\def\obeta{\overline\beta}  \def\bbA{{\mathbb A}} \def\bbB{{\mathbb B}}
 \def\bbD{{\mathbb D}} \def\bbI{{\mathbb I}} \def\bbJ{{\mathbb J}} \def\bbK{{\mathbb K}}
  \def\bbM{{\mathbb M}} \def\bbR{{\mathbb R}} \def\bbT{{\mathbb T}}
 \def\bbX{{\mathbb X}}  \def\obB{{\overline\bbB}} \def\obD{{\overline\bbD}}
 \def\obR{{\overline\bbR}}  \def\obX{{\overline\bbX}}
 \def\unu{{\underline\nu}}
\def\ka{{\kappa}}
\def\urho{{\underline\rho}} \def\usi{{\underline\sigma}} \def\rP{{\rm P}} \def\rI{{\rm I}} \def\rp{{\rm p}} 
\def\bka{{\mbox{\boldmath$\kappa$}}} \def\bsi{{\mbox{\boldmath$\sigma$}}}
\def\bfta{{\mbox{\boldmath$\eta$}}}
\def\bpi{{\mbox{\boldmath$\pi$}}} \def\bXi{{\mbox{\boldmath$\Xi$}}}
\def\Gam{{\Gamma}} \def \ux{{\underline x}}  \def\uy{{\underline y}} 
\def\ualp{{\underline\alpha}}  \def\cA{{\mathcal A}} \def\cB{{\mathcal B}} 
\def\cK{{\mathcal K}}  \def\cZ{{\mathcal Z}} \def\ocZ{\overline\cZ} 

\def\beac{\begin{array}{c}} \def\beal{\begin{array}{l}} \def\bear{\begin{array}{r}}
\def\beacl{\begin{array}{cl}} \def\beall{\begin{array}{ll}} \def\ena{\end{array}}
\def\pt{{\partial}}

\begin{description}
\item[(a)]\label{1.0}A finite collection $\Lam$ of sites;  a first-come-first-served infinite-buffer
single-server queue at each site $j\in\Lam$. 
\item[(b)] Two vectors:  
$\ulam =(\lam_i)$, $\lam_i\geq 0$ being a Poisson arrival intensity at 
site $i\in\Lam$, and $\unu =(\mu_k)$, $\mu_k> 0$ being the service intensity at site $k\in\Lam$. 
Arrivals at different sites are independent.
\item[(c)] A sub-stochastic matrix
${\tP}=(\rp_{ik},\,i,k\in\Lam )$: 
$\rp_{ik}\geq 0,$ $\sum_{k\in\Lam}\rp_{ik}\leq 1$.
After completing service at site $i$, a task is transferred to site $k$ with probability $\rp_{ik}$  
and leaves the network with probability $\rp^*_i=1-\sum_{k\in\Lam}\rp_{ik}$.

The above description gives rise to a continuous-time MP with states 
$\unn=(n_i,\,i\in\Lam)\in\Z_+^\Lam$ where $n_i\in\Z_+:=\{0,1,2,\ldots\}$; this will be a prototype of
the basic MP ${\wt X}^{(z)}(t)$. The generator 
matrix $\bQ=(Q(\unn ,\unn^\prime))$ of the process has non-zero entries corresponding to 
the following transitions:
\beq\label{2.2}\beall
{Q}(\unn ,\unn +\ue^i)=\lam_i&\hbox{an arrival of a task at site $i$,}\\
{Q}(\unn ,\unn -\ue^i)=\mu_i\rp^*_i{\mathbf 1}(n_i\geq 1)&\hbox{a task exits 
from site $i$ out of the network,}\\
{Q}(\unn ,\unn +\ue^{i\to k})=\mu_i\rp_{ik}{\mathbf 1}(n_i\geq 1)&\hbox{a task jumps from site $i$ 
to $k$.}\ena\eeq
Here $\ue^i=(e^i_l,\,l\in\Lam)\in\Z_+^\Lam$ \ has \ $e^i_l=\delta_{il}$, $\ue^{i\to k}=\ue^k-\ue^i$, \ and  
${\mathbf 1}$ stands for the indicator.

\item[(d)] Let $\;\tt I\;$ denote the unit matrix. The total intensities $\rho_i$ of the flows through sites $i\in\Lam$ form a
vector $\urho =(\rho_i)$ where
\beq\label{2.11}\urho =\ulam ({\tt I}-{\tt P})^{-1},\;\hbox{where ${\tt I}-{\tt P}$ is supposed to be invertible.}\eeq 
\item[(e)] The IM $\nu$ for the process is given by $\nu (\unn )=\prod_{i\in\Lam}\left(\rho_i/\mu_i
\right)^{n_i}$, $\unn =(n_i)$.\end{description}
The sub-criticality condition (SCC) reads:
$\diy\rho_i/\mu_i<1$, $i\in\Lam$, and 
gives rise to an EPD $\pi$ in the form of the product of geometric marginals with parameter 
$\rho_i/\mu_i$:
\beq\label{2.14}\hbox{$\pi (\unn )=\left[\prod_{i\in\Lam}\left(\rho_i/\mu_i\right)^{n_i}\right]/\bXi,\;\;
\unn =(n_i)$, \;where $\;\bXi =\sum_{\unn\in\bbZ_+^\Lam}\pi (\unn )=
\prod_{i\in\Lam}\left(1-\rho_i/\mu_i\right)$.}\eeq
In fact, assuming that matrix $\tt P$ is irreducible (i.e., ${\tt P}^s$ has strictly positive entries
for some positive integer $s$) and that the SCC holds, the JN MP is positive recurrent, and 
$\pi$ is the unique EPD.
\vskip .5 truecm

{\bf 2.B. The combined generator.} \ We now give the description of the model modifying the 
basic example from \cite{GPSY}. Let $\cZ$ be a finite or  a countable set. 
The state of the combined MP is a pair $(z,\unn )$ where 
$z\in\cZ$ indicates a state of the environment (SE) and $\unn\in\bbZ^\Lam$ a state of a basic 
process. Furthermore, given $z\in\cZ$,
we fix a collection of triples $(\ulam^{(z)},\unu^{(z)},{\tt P}^{(z)})$, with $\ulam^{(z)}=(\lam^{(z)}_i)$, 
$\umu^{(z)}=(\mu^{(z)}_i)$, ${\tt P}^{(z)}=(\rp^{(z)}_{ik})$,  
assuming that, $\forall$ $z\in\cZ$, Eqns \eqref{2.2}--\eqref{2.11} hold, and matrix ${\tt P}^{(z)}$
is irreducible. The jump rates form an (infinite) generator matrix $\BR 
=\Big(R[(z,\unn ),(z^\prime ,\unn^\prime )]\Big)$. The entries $R[(z,\unn ),(z^\prime ,\unn^\prime )]$ are:
\beq\label{eq:2.15}\begin{array}{ll}
{R}[(z,\unn ),(z,\unn +\ue^i)]=\alpha (z)\lam^{(z)}_i&\hbox{a task arrival at site $i$ under SE $z$,}\\
{R}[(z,\unn ),(z,\unn -\ue^i)]=\alpha (z)\mu^{(z)}_i{\rp_i^*}^{(z)}{\mathbf 1}(n_i\geq 1)&\hbox{a task 
exit from site $i$ under SE $z$,}\\
{R}[(z,\unn ),(z,\unn +\ue^{i\to k})]=\alpha (z)\mu^{(z)}_i\rp^{(z)}_{ik}{\mathbf 1}(n_i\geq 1)&\hbox{a task 
moves $i\to k$ under SE $z$,}\\
{R}[(z,\unn ),(z^\prime ,\unn )]=\sigma (z)\tau^{(\uNn )} (z,z^\prime )\prod\limits_{i\in\Lam}
\left(\rho^{(z)}_i/\mu^{(z)}_i\right)^{-n_i}
&\hbox{the SE changes from $z$ to $z^\prime$,}\\
{R}[(z,\unn ),(z^\prime ,\unn^\prime)]=0&\hbox{$\forall$ other pair $(z^\prime ,\unn^\prime)
\in\cZ\times\bbZ_+^\Lam$}\\
\;&\hbox{with $(z,\unn)\neq (z^\prime ,\unn^\prime)$,}\\
{R}[(z,\unn ),(z,\unn )]=-\sum\limits_{(z^\prime ,\unn^\prime)
\in\cZ\times\bbZ_+^\Lam\setminus\{(z,\unn )\}}{R}[(z,\unn ),(z^\prime ,\unn^\prime)].
\ena \eeq
Here $(\alpha (z),\;z\in\cZ)\in\bbR_+^\cZ$ and $(\sigma (z),\;z\in\cZ)\in\bbR_+^\cZ$ are given vectors with 
entries $\alpha (z)\in [0,\infty )$, $\sigma (z)\in (0,\infty )$; they represent two forms
of time-scaling: for the basic MP under SE $z$ and for the exit from $z$. Next,
$\rT^{(\uNn )}=(\tau^{(\uNn )}(z,z^\prime ))$ is the nominal SE-jump intensity matrix (possibly, depending
on $\unn$) for which we assume that 
\beq\label{eq:2.16A}\hbox{$\tau^{(\uNn )}(z,z^\prime )\geq 0$, $\;\sum_{z^\prime\in\cZ}\tau^{(\uNn )}
(z,z^\prime )
=\sum_{z^\prime\in\cZ}\tau^{(\uNn )}(z^\prime , z)$, $\;\tau (z,z)=0$, $\;\forall\;\;z\in\cZ$, 
$\unn\in\bbZ_+^\Lam$.}\eeq
  
One can see that the top three lines in \eqref{eq:2.15} emerge from the jump rates of a basic
MP ${\wt X}^{(z)}(t)$ whereas the bottom line is related to the SE MP ${\wt Z}(t)$.    
  
We say that measure $\bfta$ on $\cZ\times\bbZ_+^\Lam$ and the collection of rates 
$\BR=\{R[(z,\unn ),(z^\prime,\unn^\prime)]\}$
satisfy the WIE (the weak invariance equation) if,  
$\forall$ $(z,\unn )\in\cZ\times\bbZ_+^\Lam$, 
\beq\label{eq:WIE1}\hbox{$\sum\limits_{(z^\prime ,\uNn^\prime )\neq (z,\uNn)}
\bfta (z,\unn )R[(z,\unn ), (z^\prime ,\unn^\prime )]
=\sum\limits_{(z^\prime ,\uNn^\prime )\neq (z,\uNn)}\bfta (z^\prime ,\unn^\prime )
R[(z^\prime ,\unn^\prime ),(z,\unn )]$.}\eeq     
  
{\bf Theorem 2.1.} (I) {\sl Set: 
\beq\label{eq:2.16B}\hbox{$\bka (z,\unn )=\left[\prod_{i\in\Lam}
\left(\rho^{(z)}_i/\mu^{(z)}_i\right)^{n_i}\right]\big/\sigma (z),\;\;z\in\cZ,\;\unn =(n_i)\in\bbZ^\Lam$.}\eeq
Then measure $\bka$ and the collection of rates from Eqn \eqref{eq:2.15} satisfy the WIE. }

(II) {\sl Assuming the SCC 
\beq\label{eq:2.17}\diy\rho^{(z)}_i/\mu^{(z)}_i<1\;\;\forall\;\;i,z\;\hbox{ and \ \
$\bXi =\sum_{z\in\cZ}
\left[\prod_{i\in\Lam}\left(1-\rho^{(z)}_i/\mu^{(z)}_i\right)\right]\big/\sigma (z)<\infty$}\eeq
yields a PM $\bpi$ satisfying the WIE: $\bpi (z,\unn )=\left[\prod_{i\in\Lam}
\left(\rho^{(z)}_i/\mu^{(z)}_i\right)^{n_i}\right]\big/[\bXi (\Lam )\sigma (z)]$.}
 
{\it Proof.} Assertion (II) follows from (I) so we focus on the proof of (I). We check  
partial balance equations: $\forall$ $(z,\unn )$, we assert that $F^{\rm{out}}(z,\unn )=F^{\rm{in}}(z,\unn )$. Here
$$\hbox{$F^{\rm{out}}(z,\unn )= \sum_{(z^\prime ,\unn^\prime )}\bka (z,\unn )R[(z,\unn ), (z^\prime ,\unn^\prime )]$,
$F^{\rm{in}}(z,\unn )=
\sum_{(z^\prime ,\unn^\prime )}\bka (z^\prime ,\unn^\prime )R[(z^\prime ,\unn^\prime ),(z,\unn )]$.}$$   
It is convenient to represent $F^{\rm{out}}(z,\unn )=F^{\rm{out}}_1(z,\unn )+F^{\rm{out}}_2(z,\unn )$, $F^{\rm{in}}(z,\unn )=F^{\rm{in}}_1(z,\unn )+F^{\rm{in}}_2(z,\unn )$, 
and prove that $F^{\rm{out}}_1(z,\unn )=F^{\rm{in}}_1(z,\unn ),\;\;F^{\rm{out}}_2(z,\unn )=F^{\rm{in}}_2(z,\unn )$, with 
$$\hbox{$F^{\rm{out}}_1(z,\unn )= \sum_{\uNn^\prime }\bka (z,\unn )R[(z,\unn ), (z,\unn^\prime )]$, $
F^{\rm{out}}_2(z,\unn )= \sum_{z^\prime }\bka (z,\unn )R[(z,\unn ), (z^\prime ,\unn )]$,}$$
and
$$\hbox{$F^{\rm{in}}_1(z,\unn )=
\sum_{\unn^\prime }\bka (z,\unn^\prime )R[(z,\unn^\prime ),(z,\unn )]$, $
F^{\rm{in}}_2(z,\unn )=
\sum_{z^\prime }\bka (z^\prime ,\unn )R[(z^\prime ,\unn),(z,\unn )]$.}$$  

After omitting the factor $1/\sigma (z)$, the equation $F^{\rm{out}}_1(z,\unn )=F^{\rm{in}}_1(z,\unn )$ 
means that 
$$\hbox{$\sum_{\unn^\prime }\nu^{(z)}(\unn )Q^{(z)}(\unn , \unn^\prime )=
\sum_{\unn^\prime }\nu^{(z)}(\unn^\prime )Q^{(z)}(\unn^\prime ,\unn )$}$$
which holds as $\nu^{(z)}$ is an IM for $Q^{(z)}$. Next, $F^{\rm{out}}_2(z,\unn )
=F^{\rm{in}}_2(z,\unn )$ is equivalent to \eqref{eq:2.16A} since
$$\hbox{$F^{\rm{out}}_2(z,\unn )=\sum_{z^\prime }\big(1/\sigma (z)\big)\prod_{i}\left(\rho^{(z)}_i
/\mu^{(z)}_i\right)^{n_i} 
\tau^{(\uNn )} (z,z^\prime )\prod_{i}\left(\rho^{(z)}_i/\mu^{(z)}_i\right)^{-n_i}\sigma (z)=
\sum_{z^\prime }\tau^{(\uNn )} (z,z^\prime )$}$$ 
and similarly $ F^{\rm{in}}_2(z,\unn )=\sum_{z^\prime }\tau^{(\uNn )}(z^\prime , z)$.  $\qquad\Box$
\vskip .5 truecm

{\bf Remarks.\ 2.1.} Note that values $\alpha (z)$ and $\tau^{(\uNn )}(z,z^\prime)$ do not enter  expression 
\eqref{eq:2.16B}. However, if $\alpha (z)\equiv 0$ then $\forall$ fixed $\unn$ pairs $(z,\unn)$, $z\in\cZ$,
form a closed communicating class supporting an IM $\bka_0=\bka_0^{(\uNn)}$ with values
$\bka_0(z,\uNn)=1/\sigma (z)$.   


{\bf 2.2.} Under additional conditions of a `moderate growthÕ of functions $\alpha (z)$ and 
$\sigma (z)$ (cf. Theorem 3.1, P. 376 and Corollary 3.2, P. 379 in \cite{EK}),  the collection of rates 
\eqref{eq:2.15} defines a combined MP $(Z(t),X(t))$ 
on $\cZ\times\bbZ_+^\Lam$ with a Feller semi-group of transition operators. (Here it means that,
$\forall$ $t>0$, the transition matrix of the process takes the space of bounded functions 
$(z,\unn)\mapsto \phi (z,\unn )$ to itself).  
Physically speaking, for `nice` $\alpha (z)$ and $\sigma (z)$  we obtain a non-explosive combined 
MP. (For a finite set $\cZ$  
this is automatically true.)

Formally, a sufficient condition (deduced from the aforementioned results in \cite{EK})
is that the following three assumptions (i)--(iii) hold. (i) For 
$\oR (z,\unn )=-R[(z,\unn ),(z,\unn )]$ we assume:
\beq\label{eq:ovR}\beal\diy\operatornamewithlimits{\sup}\limits_{(z,\unn )\in\cZ\times\bbZ_+^\Lam}
\oR (z,\unn )<\infty\;\hbox{ where }\\
\qquad\qquad\diy\oR (z,\unn )=\sigma (z)\sum\limits_{z^\prime\in\cZ}\tau^{(\uNn )}_{zz^\prime} 
\prod\limits_{i\in\Lam}\left(\frac{\rho^{(z)}_i}{\mu^{(z)}_i}\right)^{-n_i}
+\alpha (z)\sum\limits_{i\in\Lam}\big[\lam^{(z)}_i+ \mu^{(z)}_i
{\mathbf 1}(n_i\geq 1)\big].\ena\eeq
Next, (ii) after compactifying space $\cZ\times\bbZ_+^\Lam$ by a point $\Delta$, we 
suppose that $\forall$ finite subset $\cK\in\cZ\times\bbZ_+^\Lam$,  
\beq\lim\limits_{(z,\unn )\to\Delta}\sum\limits_{(z^\prime ,\unn^\prime)\in\cK}
\big|R[(z,\unn),(z^\prime ,\unn^\prime )]\big| =0,\eeq
Finally, (iii) 
we assume that  
\beq\diy\operatornamewithlimits{\sup}\limits_{(z,\unn )\in\cZ\times\bbZ_+^\Lam}
\;\sum\limits_{(z^\prime ,\unn^\prime )\in\cZ\times\bbZ_+^\Lam} \Big|\Big\{\oR (z,\unn )
-\oR (z^\prime,\unn^\prime ) 
\Big\}\,R[(z,\unn ),(z^\prime ,\unn^\prime)]\Big| <\infty .\eeq

Results about
stationary distributions (see Proposition 9.2, P. 239, from \cite{EK}) imply that under assumptions
(i)--(iii), 
a measure satisfying
the WIE is invariant under the transition semi-group.  Furthermore, under 
assumptions (i) -- (iii) and condition \eqref{eq:2.17}, 
the combined MP $(Z(t),X(t))$ 
is positive recurrent.

{\bf 2.3.} In assumption (i) (cf. Eqn \eqref{eq:ovR}), the troublesome summand is 
$$\oR^{(1)}(z,\unn):=\diy\sigma (z)\sum\limits_{z^\prime\in\cZ}\tau^{(\uNn )}_{zz^\prime} 
\prod\limits_{i\in\Lam}\left(\frac{\mu^{(z)}_i}{\rho^{(z)}_i}\right)^{n_i}.$$
A sufficient condition for $\operatornamewithlimits{\sup}\limits_{(z,\unn )\in\cZ\times\bbZ_+^\Lam}
\oR^{(1)}(z,\unn )<\infty$ covering  a host of realistic situations is that (iv) $\tau^{(\unn )}_{zz^\prime}$
is of the form $\tau^{(\unn )}_{zz^\prime}=h(\unn )\,\otau_{zz^\prime}$ where \ $\otau_{zz^\prime}$ 
does not
depend on \ $\unn$ \ and has $\operatornamewithlimits{\sup}\limits_{z\in\cZ}
\sum\limits_{z^\prime\in\cZ}\otau_{zz^\prime}<\infty$
and \ (v) $h(\unn )$ is such that $\diy\operatornamewithlimits{\sup}\limits_{z\in\cZ,\uNn\in\bbZ_+^\Lam}
h(\unn )\;\prod_{i\in\Lam}\left(\frac{\mu^{(z)}_i}{\rho^{(z)}_i}\right)^{n_i}<\infty$. For example, 
assuming that 
$S(n):=\diy\operatornamewithlimits{\sup}\limits_{z\in\cZ}\prod\limits_{i\in\Lam}
\;\left(\frac{\mu^{(z)}_i}{\rho^{(z)}_i}\right)^{n_i}<\infty$, the value $h(\unn)$ can be 
selected as $S(n)^{-1}$.

\section{A simple exclusion model in a random environment} \label{sec:3}

In this section, we focus on a class of MPs with local interaction arising from a symmetric
simple exclusion model; cf.\cite{Lig1}--\cite{Lig3}. The state $x$ of the basic MP is 
considered as a spin-$0,1$ configuration (aka an occupancy configuration)  on a cubic lattice 
$\bbZ^d$, i.e., as a function 
$i\in\bbZ^d\mapsto\rx_i\in\{0,1\}$. The value $\rx_i=1$ means that site $i$ is occupied by a `particle`
whereas $\rx_i=0$ means that site $i$ is vacant/empty. (We also will use the term a `light particle` as
opposite to a `heavy` particle introduced below.) Let us fix parameters 
\beq\beac\vphi\in\bbR, \;\lam ,\mu>0,\;
\beta_{i,i^\prime}=\beta_{i^\prime i}\qquad\qquad\qquad\qquad{}\\
\qquad\qquad\qquad\hbox{ where }i,i^\prime\in\bbZ^d,\; \beta_{ii}=0\;\hbox{ and }\;
\operatornamewithlimits{\sup}\limits_{i\in\bbZ^d}\sum\limits_{i^\prime\in\bbZ^d}
\beta_{ii^\prime}<\infty .\ena\eeq
Next, we take $\bbX =\{0,1\}^{\bbZ^d}$ with the product topology (and -- when necessary --
with a metric generating this topology), 
and set $\cZ=\bbT$ where $\bbT\subset\bbZ^d$ is a finite set. Further, for 
$x=(\rx_i,\,i\in\bbZ^d)\in\{0,1\}^{\bbZ^d}$
and $z\in\bbT$ we write:  
\beq\label{eq:LigL}\beac\tL^{(z)}g(x)=\sum\limits_{i,i^\prime\in\bbZ^d:\,i\neq i^\prime\neq z}
\rx_1(1-\rx_{i^\prime})\beta_{i,i^\prime}\big[g(x+e^{i\to i^\prime})-g(x)\big]
\qquad\qquad\qquad{}|\\
+\sum\limits_{i\in\bbZ^d:\,i\neq z}
\Big\{\rx_i(1-\rx_i)\theta_{iz}e^\vphi\big[g(x+e^{i\to z})-g(x)\big]\\
\qquad\qquad +\lam\big[g(x+e^{i})-g(x)\big]
+\rx_i\mu\big[g(x-e^{i})-g(x)\big]\Big\}\\
\qquad +\lam e^\vphi \big[g(x+e^{z})-g(x)\big]
+\rx_z\mu\big[g(x-e^{z})-g(x)\big].\ena\eeq
Here $e^{i}\in\{0,1\}^{\bbZ^d}$ is the configuration 
where all values of spin are $0$ except for that at site $i$, which is $1$, and we set: 
$e^{i\to i^\prime}:=e^{i}
-e^{i^\prime }$. (This covers the case where $i$ or $i^\prime$ coincides with $z$.) Adding and 
subtracting configurations means addition and subtraction of functions. 

A construction below repeats that from \cite{Lig1}, Sect. I.3. Denote by $\rC_\rL(\{0,1\}^{\bbZ^d})$ the 
space of continuous functions $\rC(\{0,1\}\times\bbZ^d)$. 
Operator $\tL^{(z)}$ in \eqref{eq:LigL} acts initially on the space $\rC_\rL(\{0,1\}^{\bbZ^d})$ of 
Lipschitz-type functions $x\in\{0,1\}^{\bbZ^d}\mapsto g(x)$:
\beq\label{eq:CLz}\bear\rC_\rL(\{0,1\}^{\bbZ^d})=\bigg\{g\in\rC (\{0,1\}^{\bbZ^d}):\;|||g|||
:=\operatornamewithlimits{\sup}\limits_{x\in\{0,1\}^{\bbZ^d}}\bigg[
\sum\limits_{i\in\bbZ^d}\rx_i\left|g\left(x-e^i\right)-g(x)\right|\quad\qquad{}\\
+\sum\limits_{i\in\bbZ^d}(1-\rx_i)\left|g\left(x+e^i\right)-g(x)\right|
+\sum\limits_{i,i^\prime\in\bbZ^d}\rx_i(1-\rx_{i^\prime})
\left|g\left(x+e^{i\to i^\prime}\right)-g(x)\right|\bigg]<\infty\bigg\}.\ena\eeq
Next, we extend \eqref{eq:LigL} to a closed operator in 
$\rC (\{0,1\}^{\bbZ^d})$. 
The closed operator is still denoted by $\tL^{(z)}$; for its domain 
we use the notation $\rD (\tL^{(z)})$. By construction, $\rC_\rL(\{0,1\}^{\bbZ^d})$ is a core for $\tL^{(z)}$.
Physically, $\tL^{(z)}$ is a generator of an MP ${\wt X}^{(z)}$ on the state space $\{0,1\}^{\bbZ^d}$ representing
an `open` simple-exclusion model, in presence of a (single) `heavy particle` placed at site $z\in\bbT$. In 
this model, 
a light particle can jump from site $i$ to site $i^\prime$ at rates $\beta_{ii^\prime}$ or
$\beta_{ii^\prime}e^\vphi$, depending on the status of sites $i$ and $i^\prime$ (vacant, occupied by a light 
particle, occupied/not occupied by a heavy particle). In any case, jumps are performed only if
the simple-exclusion restriction is respected: at  most one light particle at any given site. (A simultaneous
presence of a light and the heavy particle at a given site is allowed.)
A light particle can also 
be annihilated at rate $\mu$ and created at rates $\lam$ or $\lam e^\phi$, depending on whetger
the given site is occupied by the heavy particle. Factor $e^\phi$ indicates an impact
that a heavy particle has on the dynamics of a light-particle configuration $x$: when $\vphi >0$, the 
heavy particle attracts the light ones, when $\vphi <0$, it repels them. Following \cite{Lig1}--\cite{Lig3}, 
it is possible to check that $\forall$ $z\in\bbT$ there exists a unique Feller semi-group of operators in   
$\rC (\{0,1\}^{\bbZ^d})$ generated by $\tL^{(z)}$. See, e.g., \cite{Lig1}, Theorem I.3.9, P. 27,
or  \cite{Lig3}, Theorem 4.68. This yields a basic MP ${\wt X}^{(z)}$ in $\{0,1\}^{\bbZ^d}$.

Furthermore, an IM $\nu^{(z)}$ (in fact, an EPD) for process ${\wt X}^{(z)}$ is given by
\beq\label{eq:IMSE}
\nu^{(z)}=\left(\prod\limits_{i\in\bbZ^d:\,i\neq z}\tP^{(i)}\right)\times\tQ^{(z)}.\eeq
It is a product-measure (aka an inhomogeneous Bernoulli measure) where the spin values are independent 
for different sites, and the marginal
distribution for an individual spin is either $\tP^{(i)}\simeq\tP$ or $\tQ^{(i)}\simeq\tQ$,
depending on whether $i\neq z$ (i.e., site $i$ not occupied by a heavy particle) or $i=z$ (i.e., $i$ 
contains the heavy particle). Both $\tP$ and $\tQ$ are probability distributions on the two-outcome 
set $\{0,1\}$:
\beq\label{eq:tPtQ}\tP(1)=\frac{\lam}{\lam+\mu},\;\tP(0)=\frac{\mu}{\lam+\mu},\;\tQ(1)
=\frac{\lam e^\vphi}{\lam e^\vphi+\mu},\;\;\tQ(0)=\frac{\mu}{\lam e^\vphi +\mu}.\eeq
All measures $\nu^{(z)}$ are absolutely continuous relative to 
$\gam =\operatornamewithlimits{\times}\limits_{i\in\bbZ^d}\tP^{(i)}$, with the Radon--Nikodym densities 
\beq\label{eq:mzxSE}\bear\diy m(z,x)=\frac{\rd\nu^{(z)}(x)}{\rd\gam (x)}=
\frac{\lam +\mu}{\lam e^\vphi +\mu}\Big[e^\vphi{\mathbf 1}(\rx_z=1)+{\mathbf 1}(\rx_z=0)\Big]
=\frac{\lam +\mu}{\lam e^\vphi +\mu}e^{\rx_z\phi}\qquad\qquad{}\\
\hbox{where }z\in\bbT\;\hbox{ and }\;x =(\rx_i)\in\{0,1\}^{\bbZ^d}.\ena\eeq
We also have that $\int_{\{0,1\}^{\bbZ^d}}\tL^{(z)}g(x)\rd\nu^{(z)}(x)=0$, $\forall$ $z\in\bbT$ and 
$g\in\rC_\rL(\{0,1\}^{\bbZ^d})$.

Next, the SE (state of environment) process ${\wt Z}(t)$ has the generator $\tA$ defined by a finite matrix for
a function/vector $z\in\bbT\mapsto f(z)$,
\beq\label{eq:LigA}\bear\tA f(z)=\sum\limits_{z^\prime\in\bbT}\tau_{zz^\prime}
[f(z^\prime )-f(z)]\hbox{ where rates }\;\tau_{zz^\prime}\;\hbox{ obey}\qquad\qquad{}\\
\tau_{zz^\prime}=\tau_{z^\prime z}\geq 0,\;\hbox{ and }\;\tau_{zz}=0.
\ena\eeq
Process ${\wt Z}(t)$ is  a random walk in $\bbT$ with a counting invariant measure
$\gam$: $\gam (\bbD)=\#\,\bbD$, \ $\bbD\subseteq\bbT$. The invariance property simply means that 
$\sum\limits_{z^\prime}\tau_{zz^\prime}=\sum\limits_{z^\prime}\tau_{z^\prime z}$ and is deduced
from the symmetry condition $\tau_{zz^\prime}=\tau_{z^\prime z}$ in
\eqref{eq:LigA}. It implies that 
$\sum\limits_{z\in\bbT}\tA f(z)=0$ for any function $f$. 

Further, consider a measure $\bka$ on $\bbT\times\{0,1\}^{\bbZ^d}$:
\beq\label{eq:Ligbka}\bka (z,\bbA )=\frac{1}{\sigma (z)}\int_{\bbA}e^{\rx_z\vphi}\rd\gam (x),\;\;
z\in\bbT,\;\;\bbA\subseteq\{0,1\}^{\bbZ^d}.\eeq

The combined generator $\BR$ is constructed from the action 
\beq\label{eq:LigM}\BR\phi (z,x)=\alpha (z)\BL^{(z)}\phi (z,x)+\frac{\sigma (z)}{e^{\rx_z\vphi}}
\BA \phi (z,x).\eeq
Here $\alpha (z)>0$ and $\sigma (z)>0$ are time-scaling coefficients.  Next, $\BL^{(z)}$ acts
on the section map $x\in\{0,1\}^{\bbZ^d}\mapsto\phi (z,x)$:
\beq\label{eq:LigBL}\beac\BL^{(z)}\phi (z,x)=\sum\limits_{i,i^\prime\in\bbZ^d:\,i\neq i^\prime\neq z}
\rx_i(1-\rx_{i^\prime})\beta_{i,i^\prime}\Big[\phi (z,x+e^{i\to i^\prime})-\phi (z,x)\Big]
\qquad\qquad\qquad{}\\
\qquad\qquad +\sum\limits_{i\in\bbZ^d:\;i\neq z}\Big\{ 
\rx_i(1-\rx_z)\beta_{iz}e^\vphi
\big[\phi (z,x+e^{i\to z})-\phi (z,x)\big]\Big\}\\
\qquad\qquad\qquad\qquad +
\lam\big[\phi (z,x+e^{i})-\phi (z,x)\big]
+\rx_i\mu\big[\phi (z,x-e^{i})-\phi (z,x)\big]\Big\}\\
\qquad\qquad\qquad +\lam e^\vphi \big[\phi (z,x+e^{z})-\phi (z,x)\big]
+\rx_z\mu\big[\phi (z,x-e^{z})-\phi (z,x)\big].\ena\eeq
Further, 
\beq\label{eq:LigBA}\BA \phi (z,x)=\sum\limits_{z^\prime\in\bbZ^d}\tau_{zz^\prime}
[\phi (z^\prime ,x)-\phi (z,x)].\eeq
(The constant factor $\diy\frac{\lam e^\vphi +\mu}{\lam +\mu}$ has been absorbed into $\sigma (z)$.)

Let $\rC (\bbT\times\{0,1\}^{\bbZ^d})=\rC (\{0,1\}^{\bbZ^d})^{\bbT}$ be the space of continuous functions 
$\rC_{\rL}(\bbT\times \{0,1\}^{\bbZ^d})$ on $\bbT\times\{0,1\}^{\bbZ^d}$. As in the case of $\tL^{(x)}$, 
we consider $\BR$ initially on space $\rC_{\rL}(\bbT\times \{0,1\}^{\bbZ^d})$ of Lipschitz-type functions:  
\beq\beal\rC_{\rL}(\bbT\times \{0,1\}^{\bbZ^d})=\bigg\{\phi\in\rC(\{0,1\}^{\bbZ^d}):\;\;
\hbox{$\forall$ $z\in\bbT$,}\\  
\qquad\hbox{the section map $x\in\{0,1\}^{\bbZ^d}\mapsto\phi (z,x)$
lies in $\rC_\rL(\{0,1\}^{\bbZ^d})$; cf. \eqref{eq:CLz}}
\bigg\}.\ena\eeq
Then take the closure in $\rC (\bbT\times\{0,1\}^{\bbZ^d})$, keeping 
for the obtained closed operator the same notation $\BR$. 
By construction, $\rC_{\rL}(\bbT\times \{0,1\}^{\bbZ^d})$ is a core for $\BR$. Owing to Theorem 3.9 from
 \cite{Lig1} and/or Theorem 4.68 from \cite{Lig3}, there exists a unique 
Feller semi-group on $\rC (\bbT\times\{0,1\}^{\bbZ^d})$
generated by $\BR$; the corresponding MP is denoted by $(Z(t),X(t))$. The trajectories of process 
$(Z(t),X(t))$ are right-continuous maps $[0,\infty )\mapsto \bbT\times\{0,1\}^{\bbZ^d}$ with left limits. 
\vskip 3 truemm 

{\bf Theorem 3.1.} {\sl Under the above conditions, $\bka$ is a finite IM for process $(Z(t),X(t))$.}  
\vskip 3 truemm 

{\it Proof.} By virtue of Proposition I.6.10 in \cite{Lig1}, P. 52, it suffices to check the equation\\ $\int_{\bbT\times
\{0,1\}^{\bbZ^d}}(\BR\phi )(z,x)\rd\bka (z,x)=0$ for $\phi\in\rC_{\rL}(\bbT\times \{0,1\}^{\bbZ^d})$. But for any 
such function $\phi$, 
$$\beal\int_{\bbT\times
\{0,1\}^{\bbZ^d}}\;(\BR\phi )(z,x)\rd\bka (z,x)\\
\qquad =\sum\limits_{z\in\bbT}{\diy\frac{\alpha (z)}{\sigma (z)}}
\int_{\{0,1\}^{\bbZ^d}}(\BL^{(z)}\phi )(z,x)e^{\rx_z\vphi}\rd\gam (x)+\int_{\{0,1\}^{\bbZ^d}}
\sum\limits_{z\in\bbT}(\BA \phi )(z,x)\rd\gam (x).\ena$$
Observe that $\forall$ $z\in\bbT$, owing to the fact that the section map of $\phi$ belongs to 
$\rC_\rL(\{0,1\}^{\bbZ^d})$, 
$$\int_{\{0,1\}^{\bbZ^d}}(\BL^{(z)}\phi )(z,x)e^{\rx_z\vphi}\rd\gam (x)=\frac{\lam +\mu}{\lam e^\vphi +\mu}
\int_{\{0,1\}^{\bbZ^d}}(\BL^{(z)}\phi )(z,x)\rd\nu^{(z)}(x)=0.$$
Also, for $\gam$-a.a.  $x\in\{0,1\}^{\bbZ^d}$, 
we have that $\sum\limits_{z\in\bbT}
(\BA \phi )(z,x)=0$. Hence, $\int_{\bbT\times\{0,1\}^{\bbZ^d}}(\BR\phi )\rd\bka =0$.


\section{A general construction}\label{sec:4}

{\bf 4.A. Jump Markov processes.} \ Theorem 1 admits an extention where $(\cZ ,\ups )$ is a 
(standard) measure space. (Usual measurability and countable additivity assumptions must be adopted
throughout this sub-section.)  Here, sums $\sum_z$ and $\sum_{z^\prime}$
are replaced by integrals against measure $\ups$; measures on space $\cZ\times\bbZ^\Lam$ 
are supposed to be given via (Radon--Nikodym) densities wrt $\ups$ times the counting measure. 
Viz. in \eqref{eq:2.17} one reads $\;\bXi =\int_{\cZ}\left\{\left[
\prod_{i\in\Lam}\left(1-\rho^{(z)}_i/\mu^{(z)}_i\right)\right]\big/\sigma (z)\right\}\rd\ups (z)<\infty$. 
Existence of the combined MP $(Z(t),X(t))$ should be analysed separately; cf. Remarcs 2B, 2C.

Next, $\bbZ^\Lam$ can also be replaced with a measure space, 
$(\bbX ,\gam)$; here  the combined state space is again the Cartesian product 
$(\cZ\times\bbX,\ups\times\gam )$. 
We refer the reader to in \cite{EK}, Section 2 of Chapter 4, PP. 162--173, and Section 3 
of Chapter  8, PP. 376--382 , for a detailed treatment of 
jump MPs and their generators. Assume that for $x\in\bbX$ we have a kernel $\tau^{(x)}(z,z^\prime )\geq 0$,
$z,z^\prime\in\cZ$, and for $z\in\cZ$ a kernel $Q^{( z )}( x , x^\prime)\geq 0$ and a
function $m (z,x)>0$,  $x, x^\prime\in\bbX$, such that  $\tau^{(x)}(z,z)=Q^{(z)}(x,x)=0$ and
\beq\label{eq:3.103}\begin{array}{c}
\hbox{$\int_\cZ \tau^{(x)}(z, z^\prime )\rd\ups (z^\prime )=\int_\cZ \tau^{(x)}(z^\prime ,z)
\rd\ups (z^\prime )$,}\\ 
\hbox{$m(z,x)\int_{\bbX} Q^{(z)}(x,y)\rd\gam (y)
=\int_{\bbX}m(z,y)Q^{(z)}(y,x)\rd\gam (y)$,}\ena  \hbox{ for $(\ups\times\gam )$-a.a 
$(z,x)\in\cZ\times\bbX$.}\eeq
Physically speaking, Eqn \eqref{eq:3.103} yields that the environment MPs have an IM $\ups$ whereas 
the basic MP under SE (state of environment) $z$  has an IM $\nu^{(z)}$ with Radon--Nikodym density 
$\diy\frac{\nu^{(z)}(\rd x)}{\gam (\rd x)}=m(z,x)$. (As above, a formal 
construction of these MPs requires additional assumptions discussed in \cite{EK}.)
 
The combined generator has rates
$R[( z, x),( z^\prime, x^\prime)]$:
\beq\label{3.102}
{R}[( z , x ),( z , x^\prime)]=\alpha ( z )Q^{( z )}( x , x^\prime ),\;\;\;
{R}[( z , x ),( z^\prime, x )]=\sigma ( z )\frac{\tau^{(x)}(z ,  z^\prime)}{m(z,x )},\eeq
with time-scale coefficients $\alpha ( z )\geq 0$, $\sigma ( z )>0$ having the same meaning as 
before. We say that a density $\bfta (z,x)$ on $(\cZ\times\bbX ,\ups\times\gam)$ and kernel 
$\BR =\{R[(z,x),(z^\prime ,x^\prime )]\}$ satisfy the WIE if, for $(\ups\times\gam )$-a.a 
$(z,x)\in\cZ\times\bbX$, 
\beq\label{eq:WIE2}\begin{array}{l}\hbox{$\bfta (z, x )\int_{\cZ\times\bbX} R[(z,x)(z^\prime, x^\prime)]
{\mathbf 1}((z,x)\neq(z^\prime,x^\prime))\rd\ups (z^\prime )\rd\gam (x^\prime)$}\\
\qquad\hbox{$=\int_{\cZ\times\bbX}\bfta (z^\prime ,x^\prime ))R[(z^\prime ,x^\prime ),(z,x)]
{\mathbf 1}((z,x)\neq(z^\prime,x^\prime))\rd\ups (z^\prime )\rd\gam (x^\prime )$.}\ena \eeq

{\bf Theorem 4.1.} (I) {\sl  
The density $\bka ( z , x)=\nu^{(z)}( x )/\sigma ( z )$ considered relative to measure 
$\ups\times\gam$  and the kernel $R[( z, x),( z^\prime, x^\prime)]$ from Eqn \eqref{3.102} satisfy the 
WIE on $\cZ\times\bbX$.}
(II) {\sl Under the SCC $\;\bXi =\int_{\cZ\times\bbX}\bka (z,x)\rd\ups (z)\rd\gam (x)<\infty $, 
the PD $\bpi ( z , x )=\bka ( z , x )/\bXi$ also satisfies the WIE.}
\vskip .2 truecm

{\it Proof}: Repeats that of Theorem 2.1 by replacing sums with integrals. $\qquad\Box$

{\bf Remarks.} {\bf 4.1.} \ As before, the kernels $\tau^{(x)}(z,z^\prime )$ enter the density $\bka$ 
indirectly, through measure $\ups$ and condition \eqref {eq:3.103}.

{\bf 4.2.} \ Under additional assumptions on spaces $(\cZ ,\ups )$ and $(\bbX ,\gam )$, functions 
$\alpha (z)$ and $\sigma (z)$ and kernels 
$\tau^{(x)}(z,z^\prime )$ and $Q^{(z)}(x,x^\prime )$ (see Theorem 3.1, P. 377 in \cite{EK}), operator $\BR$ 
generates a strongly continuous Feller semi-group on $\cZ\times\bbX$. 
It can 
also be achieved that $\bka$ yields an IM for this corresponding MP.  
\vskip .3 truecm

{\bf 4.B. General assumptions on component MPs.} In sub-sections 4.B -- 4.D  we 
develop further a general construction of a combined generator.  Until the end of Section \ref{sec:4},
the spaces $\cZ$ and $\bbX$ below are supposed to be locally compact Polish, and all measures under
consideration are Borel, countably-additive and finite on compact sets. 
 
Let us proceed with formal definitions. As before, an SE is represented by a point $z\in\cZ$.
We also assume that we have a measure $\ups$ on $\cZ$ which will serve as an IM for environment 
MPs. 
  
Next, there is a space $\bbX$ given; points $x\in\bbX$ are treated as 
states of a family of MPs indexed by $z\in\cZ$ -- we again call them basic MPs. Physically speaking, 
each basic MP   
has a generator and an IM indexed by $z\in\cZ$.  
More precisely, $\forall$ $z\in\cZ$, we have a measure $\nu^{(z)}$
on $\bbX$; we suppose  that each measure 
$\nu^{(z)}$ is absolutely continuous wrt a fixed measure $\gam$, with Radon--Nikodym densities 
$m (z,x)=\diy\frac{\rd\nu^{(z)}(x)}{\rd\gam (x)}>0$ for $(\ups\times\gam )$-a.a. $(z,x)\in\cZ\times\bbX$. 

Further,  $\forall$ $z\in\cZ$  
a closed linear operator $\tL^{(z)}$ is given, acting on functions $g:\,\bbX\to\bbR$ 
forming a domain $\rD(\tL^{(z)})$ dense in $\rC_\rb (\bbX )$, the space of bounded continuous
functions on $\bbX$. Further still, we suppose that  each $\tL^{(z)}$ is a generator 
of a Feller MP on $\bbX$. That is,  $\tL^{(z)}$ is a conservative operator satisfying 
the assumptions of Theorem 2.2, P. 165 in \cite{EK} (a version of the Hille--Yosida theorem).
We also assume that  $\forall$ $g$ from a core of $\tL^{(z)}$, the integral 
$\int_\bbX(\tL^{(z)}g)(x)\rd\nu^{(z)}(x)$ exists and equals $0$. (Formally, it is only the latter assumption 
that matters in the proof of Theorem 4.2 in the next sub-section.)
  

Likewise, we suppose that $\forall$ $x\in\bbX$ there is a closed linear operator $\tA^{(x)}$
acting on functions $f:\,\cZ\to\bbR$ 
forming a domain $\rD(\tA^{(x)})$ dense in $\rC_\rb(\cZ )$, the space of bounded continuous
functions on $\cZ$. In addition, we assume that  each $\tA^{(x)}$ is a generator 
of an Feller MP on $\cZ$. We also suppose that 
$\forall$ $f$ from a core of $\tA^{(x)}$, the integral $\int_\cZ(\tA^{(x)}f)(z)\rd\ups (z)$ exists and is equal to $0$. 
(Note that measure $\ups$ serves all operators $\tA^{(x)}$.)

\vskip .5 truecm

{\bf 4.C. A combined generator.} We are now going to introduce generator $\BR$ and
analyse the related WIE.  Fix a time-scale function $z\in\cZ\mapsto\sigma (z)$ with $\sigma (z)>0$ for
$\ups$-a.a. $z\in\cZ$. Consider the following measure $\bka$ on $\cZ\times\bbX$:  
\beq\label{eq:genbka}\bka (B)=\int\frac{{\mathbf 1}((z,x)\in B)}{\sigma (z)}m(z,x)\rd\gam (x)\rd\ups (z),\;\forall\;
\hbox{ Borel }\;B\subseteq \cZ\times\bbX.\eeq
Then introduce a linear map $\BR$ acting on functions $(z,x)\in (\cZ\times\bbX)\mapsto\phi (z,x)\in\bbR$, 
$\phi\in\rD(\BR)$, where 
\beq\label{eq:genBM}\BR\phi (z,x)=\BA^{(x)}\phi (z,x)+\BL^{(z)}\phi (z,x).\eeq
Here, it is understood that $\BA^{(x)}\phi (z,x)$ results in the action of $\tA^{(x)}$ in SE-variable $z$ 
succeeded by multiplication by $1/m(z,x)$,  and $\BL^{(z)}\phi (z,x)$ 
results in $\tL^{(z)}$ acting in variable $x\in\bbX$.  Formally, for $\phi =f\otimes g$, with $\phi (z,x)=f(z)g(x)$:
\beq\label{eq:BABL}[\BA^{(x)}(f\otimes g)](z,x)=g(x)\frac{\sigma (z)(\tA^{(x)}f)(z)}{m(z,x)}, 
\;\;\big[\BL^{(z)}(f\otimes g)\big] (z,x)=\alpha (z)f(z)(\tL^{(z)}g)(x)\eeq
where $z\in\cZ\mapsto\alpha (z)>0$ is another given time-scale  function. Referring to a Banach space
$\rC_\rb (\cZ )\otimes\rC_\rb (\bbX)$, with a chosen cross-norm, 
the action of $\BR$ is then extended by linearity and continuity to general functions $\phi$.

In this paper we do not intend to go into details of formal constructions of a closed operator 
based on map $\BR$ defined in \eqref{eq:genBM}.
Nor shall we try to analyse the conditions of the Hille--Yosida theorem for such an operator.
(Still, we will refer to this operator as $\BR$ and suppose that it exists.)  
It will be assumeed that the domain $\rD(\BR )$ contains functions 
$\phi\in\rC_\rb (\cZ )\otimes\rC_\rb (\bbX)$ such that 
(i) $\forall$ $z\in\cZ$, 
function $g_{\phi ,z}:\;x\mapsto\phi (z,x)$ belongs to $\rD(\tL^{(z)})$ 
(ii) $\forall$ $x\in\bbX$, 
function $f_{\phi ,x}:\;z\mapsto\phi (z,x)$ belongs to $\rD(\tA^{(x)})$. (We will say that functions $f_{\phi ,x}$ 
and $g_{\phi ,z}$ give section maps generated by function $\phi$.) 
Pictorially, we would like to treat $\BR$ as a generator of an MP $(Z(t),X(t))$ on $\cZ\times\bbX$ which is a superposition of two 
tendencies: one is to keep an SE value $z$ intact and evolve in component $\bbX$ in
accordance with the MP generated by $\alpha (z)\tL^{(z)}$, 
the other to keep state $x$ of the base MP but change the SE value $z$ by following the MP on $\cZ$ generated 
by $\sigma (z)\tA^{(x)}$. (As examples in this paper show, under additional assumptions such a process 
can be constructed.) 

For functions $\phi$ satisfying (i), (ii), formula \eqref{eq:genBM} becomes 
$$\BR\phi (z,x)=[\BA^{(x)}f_{\phi ,x}] (z)+[\BL^{(z)}g_{\phi ,z}] (x). $$
\vskip .5 truecm

{\bf Theorem 4.2.} {\sl The following properties {\rm{(I)--(II)}} are satisfied.
{\rm{(I)}} Measure $\bka$ and operator $\BR$ from  \eqref{eq:genbka}--\eqref{eq:BABL}  obey the 
WIE in the sense that the integral $\int_{\cZ\times\bbX}(\BR\phi)\rd\bka$ exists and equals $0$ 
$\forall$ function $\phi$ for which
$$\int_\bbX \left\{\int_\cZ\Big|\, [\tA^{(x)}f_{\phi ,x}](z)\Big|\rd\ups (z)\right\}\rd\gam (x)<\infty\;\hbox{ and }\;
\int_\cZ [\tA^{(x)}f_{\phi ,x}](z)\rd\ups (z)=0,\hbox{ for}\;\gam\hbox{-a.a. $x\in\bbX$,}$$ 
and
$$\int_\cZ\frac{\alpha(z)}{\sigma (z)}\left\{\int_\bbX 
\Big|\,\big[\tL^{(z)}g_{\phi ,z}\big](x)\Big|\rd\nu^{(z)}(x)\right\}\rd\ups (z)<\infty\;\hbox{ and }\;
\int_\bbX 
\big[\tL^{(z)}g_{\phi ,z}\big](x)\rd\nu^{(z)}(x)=0,\hbox{ for}\;\ups\hbox{-a.a. $z\in\cZ$.}$$

{\rm{(II)}}  If $\Xi 
=\bka (\cZ\times\bbX)<\infty$ then $\bpi (B)=\bka (B)/\Xi$ yields a probability distribution on $\cZ\times\bbX$,
with $\int_{\cZ\times\bbX}(\BR\phi )\rd\bpi =0$ $\forall$ $\phi$ as in the above definition.}
\vskip .5 truecm

{\it Proof.} As before, we focus on assertion (I). Here, for a function $\phi $, satisfying the conditions of the 
theorem,
$$\int_{\cZ\times\bbX}(\BR\phi)\rd\bka =\int_\bbX \left\{\int_\cZ [\tA^{(x)}f_{\phi ,x}](z)\rd\ups (z)\right\}
\rd\gam (x)+\int_\cZ\frac{\alpha(z)}{\sigma (z)}\left\{\int_\bbX 
\big[\tL^{(z)}g_{\phi ,z}\big](x)\rd\nu^{(z)}(x)\right\}\rd\ups (z) .$$
Both summands in the RHS vanish owing to the assumptions.
$\qquad\Box$
\vskip .5 truecm

{\bf Remarks.\ 4.3.} As in previous results, the time-scaling factor $\alpha$ does not enter expression 
for the IM in \eqref{eq:genbka}.

{\bf 4.4.} Operators $\tA^{(x)}$ enter Eqn \eqref{eq:genbka} via a {\it common}\;\ IM\;\ $\ups$.
This fact leads to a time-homogeneous combined MP. 

{\bf 4.5.} The condition that\ $\nu^{(z)}<<\gam$ and $m(z,x)>0$ guarantees that the change of 
environment results in a `non-disruptive` continuation of the combined  MP (should it exist). 

{\bf 4.6.} \ Like earlier parts of the exposition, additional assumptions are needed if we wish the 
combined MP to exist (some of these assumptions are hard to verify).  Viz., conditions guaranteeing 
that a sum of two generators is a generator of a continuous 
contraction semi-group are listed in \cite{EK}, Theorem 7.1, P. 37 and Corollary 7.2, P. 38. Sufficient 
conditions under which an operator generates a continuous positive contraction semigroup 
are given in \cite{EK}, Theorem 2.2,  P. 165. The fact that a Feller semi-group
(i.e.,  continuous positive contraction semi-group with a conservative generator) induces
an MP with Skorokhod-type sample paths is established in \cite{EK}, Theorem 2.7, P.169 and Corollary 2.8, 
P. 170. The equivalence of the WIE and the IM property is proved, under additional assumptions, in 
\cite{EK}, Proposition 9.2, P. 239. 
See also \cite{EK}, Theorem 9.17, P. 248. With regard of a sum of two generators, other 
relevant results are contained in \cite{EK}, Section 10 of Chapter 4, PP. 253--261.

{\bf 4.7.} We would like to point at an iterative property of our construction: once constructed, a combined MP 
may be used as a basic/environment MP to produce a combined process of a `higher` level.

{\bf 4.8.} The Markov property of basic processes can be replaced by weaker assumptions covering
a broader classes of random processes with an infinite memory. This can be a topic for future research.
\vskip .5 truecm

\section{Processes with diffusive components}\label{sec:5}

Diffusion MPs play an important role in the modern theory as they are described in comprehensible 
terms and provide a broad spectrum of interesting properties. A particular feature of diffusion processes 
is that they are generated by second-order differential operators, and existence and properties of an MP 
are expressed 
here through properties of coefficients and boundary conditions. 
In this section we comment on examples of the above construction where basic MPs or an environment 
MP or both are represented by diffusions.  

{\bf 5.A. Models with a jump basic process.} \ In the current sub-section we discuss models where the basic 
process is a jump MP whereas 
the SE (state of environment) process is a diffusion. More precisely, we deal with simple, but not entirely trivial, 
examples of an isolated $M/M/1/\infty$ queue, with $\bbX =\bbZ_+:=\{0,1,2,\ldots \}$. We employ 
the notation traditionally used in the literature in this area. (So, $\lam$ stands for the Poissonian arrival rate
and $\mu$ for the service rate; condition $\lam <\mu$ is necessary and sufficient for 
the queue to be stable.)

{\bf 5.A.1.} First, take $\cZ =(\eps,1)$ where $0<\eps <1$. Here the SE-point $z=\lam$ follows a diffusion 
process with coefficient $\beta (n)>0$ and with reflections 
at the borders $\eps$ and $1$. For the basic MP we have jump rates corresponding with $\mu =1$:
$$Q^{(\lam)}(n,n+1)=\lam, \;\;Q^{(\lam)}(n,n-1)={\mathbf 1}(n\geq 1), \;\;n\in\bbZ_+.$$ 
In other words, we fix the service rate $\mu =1$ and vary the arrival rate $\lam$ between $\eps$ and $1$
as prescribed by the above  diffusion. The IM 
$\ups$ on $(\eps ,1)$ is Lebesgue: $\rd\ups (\lam )=\rd \lam$, while on $\bbZ_+$ the IM 
$\nu^{(\lam)}$ is geometric, with 
$\nu^{(\lam)} (n)=\lam^n$. (For the sake of convenience, we omit the normalizing factor.)
The combined state space is 
$(0,1)\times\bbZ_+$, and the combined generator $\BR$ acts by 
\beq\label{eq:lam}\begin{array}{cl}\BR\phi (\lam,n)&\diy =\frac{\sigma (\lam)\beta (n)^2}{2\lam^n}
\frac{\pt^2}{\pt \lam^2}\phi (\lam,n)\\
\;&\quad +\alpha (\lam)\Big\{\lam\big[\phi (\lam,n+1)-\phi (\lam,n)\big]+
\big[\phi (\lam,n-1)-\phi (\lam,n)\big]{\mathbf 1}(n\geq 1)\Big\}.\ena \eeq
The domain $\rD(\BR)$ consists of functions $(\lam, n)\mapsto \phi (\lam,n)$ that are
$\rC^2$ in variable $\lam\in (\eps ,1)$ and satisfy the Neumann boundary condition $\diy\frac{\pt}{\pt \lam}
\phi (\eps +,n)=\frac{\pt}{\pt \lam}\phi (1-,n)=0$. Assuming that (i) $\sigma (\lam )$ is continuous and strictly
positive on $[\eps ,1]$ and (ii) $\sum\limits_{n\in\bbZ_+}\beta(n)^2/\eps^n<\infty$, it is possible to check 
that operator $\BR$ from Eqn \eqref{eq:lam} generates a Feller semi-group MP on  
$\rC_\rb([\eps, 1]\times\bbZ_+)$; cf. \cite{IMcK}, Section 2.1, Chapter 2, or 
\cite{KT}, Sections 3 and 4, Chapter 15. The IMs for the corresponding MP on $(Z(t),X(t))$ with generator 
\eqref{eq:lam} are analyzed in Theorem 5.1 below.
\vskip .3 truecm
 
 {\bf Theorem 5.1.} {\sl Suppose that conditions {\rm{(i)}} and {\rm{(ii)}} are fulfilled.}
 (I) {\sl An IM on $(0,1)\times\bbZ_+$ has the form $\bka (A,n)=\int_A\diy
\frac{\lam^n}{\sigma (\lam)}\rd \lam$, $A\subseteq (\eps ,1)$.} (II) {\sl If  $\Xi 
:=\sum\limits_{n\geq 0}\int_\eps^1
\diy\frac{\lam^n}{\sigma (\lam)}\rd \lam<\infty$ (the SCC), the normalized measure 
$\bpi (A,n)=\bka (A,n)\big/\Xi$ yields an EPD, and process $(Z(t),X(t))$ is positive recurrent
in the sense that $\forall$ $n\in\bbZ_+$ and interval $A\subset (\eps, 1)$ the mean return
time to the set $A\times\{n\}$ is finite.}
\vskip .3 truecm

{\it Proof.} As above, (II) is a technical corollary, and we focus on assertion (I). Here our task is to check that
measure $\bka$ is annihilated by the conjugate $\BR^*$. 
The corresponding calculation is straightforward: the shortest proof is to  pass to the density 
$\diy\frac{\lam^n}{\sigma (\lam)}=:\frac{\bka (\rd \lam,n)}{\rd \lam}$ and  check that 
$\sum\limits_{n\geq 0}\int_\eps^1\diy\frac{\lam^n}{\sigma (\lam)}(\BR f)(\lam ,n)\rd\lam=0$ $\forall$ 
$f\in\rD (\BR)$ for which the sum $\sum\limits_{n\geq 0}\int_\eps^1\diy\frac{\lam^n}{\sigma (\lam)}
\big|(\BR f)(\lam ,n)\big|\rd\lam<\infty$. 
$\qquad\Box$

E.g., with $\sigma (\lam)\equiv 1$ (no time change in the environment MP), we have 
$\Xi =\sum\limits_{n\geq 1}\diy\frac{1-\eps^n}{n}=\infty$, and 
assertion (II) does not provide an EPD. However, for $\sigma (\lam)=\diy\frac{1}{1-\lam}$ (an 
acceleration/chaotization for $\lam\sim 1$), 
$$\Xi =\sum_{n\geq 0}\int_\eps^1{\lam^n}(1-\lam)\rd \lam=\sum_{n\geq 1}\left(\frac{1-\eps^n}{n}
-\frac{1-\eps^{n+1}}{n+1}\right)<\infty ,$$
and the construction in assertion (II) guarantees existence of an EPD. 
\vskip .5 truecm 

{\bf 5.A.2.} Our next example is a dual of the previous one: now we fix a value for the arrival rate $\lam =1$ 
and let SE $z=\mu$ follow a  Brownian motion on 
$\cZ=(1,\infty )$, with a drift $b$ and a reflection at the leftmost point $1$. The IM $\rd\ups (\mu )$
on $(1,\infty )$ is $e^{2(\mu -1)b}\rd\mu$ (which is just Lebesgue when $b =0$). The combined generator 
$\BR$ acts on $\rC^2$-functions 
$(\mu ,n)\in (1,\infty )\times\bbZ_+\mapsto \phi (\mu ,n)$ with $\diy\frac{\pt}{\pt\mu}\phi (1+,n)=0$:
\beq\label{eq:mue}\begin{array}{cl}\BR\phi (\mu ,n)&\diy =\mu^n\sigma (\mu)\left[\frac{1}{2}
\frac{\pt^2}{\pt \mu^2}\phi (\mu,n)+b\frac{\pt}{\pt \mu}\phi (\mu,n)\right]\\
\;&\quad +\alpha (\mu)\Big\{\big[\phi (\mu,n+1)-\phi (\mu,n)\big]+\mu
\big[\phi (\mu,n-1)-\phi (\mu,n)\big]{\mathbf 1}(n\geq 1)\Big\}.\ena \eeq 
In this setting, we leave open the question of existence of the combined MP, focusing instead on 
the weak invariance equation (WIE). We say that a function $(\mu, n)\in (1,\infty )\times\bbZ_+\mapsto
\eta (\mu ,n)$ (with real values $\eta (\mu ,n)\in\bbR$) satisfies the
WIE with generator $\BR$ from \eqref{eq:mue} if, $\forall$ $(\mu ,n)\in (1,\infty )\times\bbZ_+$,
the following properties hold: (i)
$$\frac{1}{2}\frac{\pt^2}{\pt\mu^2}\Big[\mu^n\sigma (\mu )\eta (\mu ,n)\Big]
-b\frac{\pt}{\pt\mu}\Big[\mu^n\sigma (\mu )\eta (\mu ,n)\Big]=0,$$
and (ii) 
$$\eta (\mu ,n)\big[1+\mu {\mathbf 1}(n\geq 1)\big]=
\eta (\mu ,n-1){\mathbf 1}(n\geq 1)+\mu\eta (\mu ,n+1).$$

{\bf Theorem 5.2.} (I) {\sl The function $(\mu, n)\in (1,\infty )\times\bbZ_+\mapsto
\ka (\mu ,n)$ of the form $\ka (\mu ,n)=\diy
\frac{e^{2(\mu-1)b}}{\mu^n\sigma (\mu)}$ satisfies the
WIE with generator $\BR$ from \eqref{eq:mue}.} (II) {\sl Under the SCC
$\Xi :=\sum_{n\geq 0}\int_1^{\infty}\diy
\frac{e^{2(\mu-1)b}}{\mu^n\sigma (\mu)}\rd \mu<\infty$, $\ka$ yields a probability density function.} 

The proof is done by a direct substitution and is omitted.

Again, with $\sigma (\mu)\equiv 1$, the partition function $\Xi =\int_1^\infty\diy\frac{\mu e^{2(\mu-1)b}}{
\mu-1}\rd\mu=\infty$, regardless of the value $b$. To obtain a PDF in assertion (II), we need to introduce
a time-scale $\sigma (\mu)$ growing at $\mu\sim 1$ and -- when $b\geq 0$ (i.e., no the drift towards $1$)
-- at  $\mu\sim \infty$.
\vskip .5 truecm

{\bf 5.A.3.} Now consider the case where the SE is a pair $(\lam, \mu)$ varies according to a joint diffusion 
in a space $\cZ$ identified as a $\pi/4$-angle
$\cA=\{(\lam ,\mu):\,\mu>\lam >0\}$, with drift $\left(\begin{array}{c}\theta\\ \theta\ena \right)$, 
reflected at the 
sides $\lam =0$ and $\lam =\mu$, along the inward normal directions. It means that generator $\tA^{(x)}=\tA$ 
acts by
\beq\label{eq:A413}\tA f(\lam ,\mu) =
\frac{1}{2}\frac{\pt^2}{\pt \lam^2}f(\lam,\mu)+\frac{1}{2}
\frac{\pt^2}{\pt \mu^2}f (\lam, \mu)
+\theta\frac{\pt}{\pt \lam}f (\lam,\mu)+\theta\frac{\pt}{\pt \mu}f (\lam, \mu).\eeq
The domain of $\tA$ consists of $\rC^2$-functions $(\lam ,\mu)\in\cA\mapsto f(\lam ,\mu)$ satisfying
the boundary conditions 
$\diy\frac{\pt}{\pt\lambda} f (0+,\mu)=0$ and $\diy\left(\frac{\pt}{\pt\lam}
-\frac{\pt}{\pt\mu}\right)f (\lam ,\mu)\Big|_{\lambda =\mu}=0$.
Here the IM $\ups$ has the Lebesgue density 
\beq\label{eq:lamuth}\frac{\rd\ups (\lam,\mu)}{\rd\lam\rd\mu}
=\exp\,[2\theta (\lam +\mu)],\;\;(\lam ,\mu)\in\cA.\eeq
(A finite IM arises iff $\theta <0$.)


Accordingly, the generator $\BR$ has the form 
\beq\label{eq:lamue}\begin{array}{l}\BR\phi (\lam ,\mu,n)\\
\diy\quad =\frac{\mu^n\sigma (\lam ,\mu)}{\lam^n}
\left[\frac{1}{2}\frac{\pt^2}{\pt \lam^2}\phi (\lam,\mu,n)+\frac{1}{2}
\frac{\pt^2}{\pt \mu^2}\phi (\lam, \mu,n)
+\theta\frac{\pt}{\pt \lam}\phi (\lam,\mu,n)+\theta
\frac{\pt}{\pt \mu}\phi (\lam, \mu,n)\right]\\
\qquad +\alpha (\lam ,\mu)\Big\{\big[\phi (\lam,\mu,n+1)-\phi (\lam,\mu,n)\big]+\mu
\big[\phi (\lam,\mu,n-1)-\phi (\lam,\mu,n)\big]{\mathbf 1}(n\geq 1)\Big\}\ena \eeq 
and acts on functions $(\lam ,\mu,n)\in\cA\times\bbZ_+\mapsto\phi (\lam ,\mu,n)$ which 
are $\rC^2$ and satisfy the above boundary conditions in variables $\lam$, $\mu$.
Repeating the previous defintion {\it mutatis mutandis} yields
\vskip .3 truecm

{\bf Theorem 5.3.} (I) {\sl Consider the function $(\lam ,\mu ,n)\in\cA\times\bbZ_+\mapsto\ka (\lam ,\mu ,n)$ 
of the form 
$$\ka (\lam \mu ,n)=\diy
\frac{\lam^n\exp\,[2\theta (\lam +\mu)]}{\mu^n\sigma (\lam ,\mu)}.$$ 
Then $\ka$ satisfies the WIE with the generator $\BR$ from \eqref{eq:lamue}.}

(II) {\sl Assuming the SCC \
$\Xi :=\sum_{n\geq 0}\int_\cZ\diy
\frac{\lam^n\exp\,[2\theta (\lam +\mu)]}{\mu^n\sigma (\lam ,\mu)}\rd\lam\rd\mu<\infty$, the normalized 
function $\ka (\lam ,\mu ,,n)\big/\Xi$ yields a PDF on $\cA\times\bbZ_+$.}
\vskip .5 truecm

{\bf Remarks.\ 5.1.} The indicated form of IM $\ups$ on $\cZ=\cA$ in Eqn \eqref{eq:lamuth}
can be obtained by observing that 
the SE diffusion process is a projection, upon $\cA$, of a `covering` MP living in a quadrant 
$\cB=\{(\lam ,\mu):\lam ,\mu>0\}$. (The projection is $(\lam ,\mu)\in\cB\mapsto (\lam\wedge\mu,
\lam\vee\mu)$ where $\lam\wedge\mu=\min [\lam,\mu]$, $\lam\vee\mu=\max [\lam,\mu]$.) The 
covering MP is a diffusion with the same drift vector (parallel to
the bissectrice) and the normal reflections from the sides $\lam =0$ and $\mu=0$. It is easy to see
that the covering diffusion is a product of two 1D diffusions, one in $\lam$ and the other in $\mu$,
each on $\bbR_+=(0,\infty )$, with drift $\theta$, reflection at the origin. A more general class of 
boundary conditions can also be considered, by following results from \cite{H}--\cite{HW}, \cite{VW}.    

{\bf 5.2.} The methodology developed thus far in this section allows us to proceed with the case of a 
Jackson network. Here it is convenient to pass to vectors $\urho$ and $\umu$ and consider 
(individual or joint) diffusions in the corresponding wedge-like domains. (Matrix $\tt P$ can 
also be varied in its own simplex-type domain.) The resulting picture essentially looks like the one 
above; detailed questions are left for/can be a subject of future research.
\vskip .5 truecm

{\bf 5.B. Models with a basic diffusion process.} Models with basic diffusion MPs and a jump change of 
the environment are quite popular in some chapters of the theory of controlled diffusion processes 
(switching diffusions); see \cite{ABG}, \cite{MY}, \cite{YZ}. Here
we discuss a straightforward example where $\bbX =\bbR$. The basic MP lives in a line and has diffusion 
coefficient $1$ and a drift $z\in\bbR$; the latter is considered as an SE. Consequently, 
$\rd\ups (x)=\rd x$, and the IM $\nu^{(z)}(\rd x)=e^{2zx}\rd x$ with $m(z,x)=e^{2zx}$, $x\in\bbR$.
 Next, we take $\cZ=\bbR$ and consider a jump environment MP with a jump measure 
 $T^{(x)}(z,\rd z^\prime)$.
Let $\ups$ be an IM with $\int_\bbR\ups (\rd z)T^{(x)}(z,\rd z^\prime ) =\ups (\rd z^\prime )$.  The individual 
components become 
\beq\label{eq:L431}\tL^{(z)}g(x)=\frac{1}{2}\frac{\rd^2}{\rd x^2}g(x)+z\frac{\rd}{\rd x}g(x),\;\;
m(z,x)=e^{2zx},\;\;z,x\in\bbR , \eeq
\beq\label{eq:A431}\tA^{(x)}f(z)=\int_\bbR T^{(x)}(z,\rd z^\prime )[\phi (x,z^\prime)-\phi (x,z)],\;\;
z,x\in\bbR .\eeq
Then set
\beq\label{eq:431}\begin{array}{c}\rd\bka (z,x)= e^{2vz}\ups (\rd z)\rd x,\\
\BR\phi (z,x)=\alpha (z)\left[{\diy\frac{1}{2}\frac{\pt^2}{\pt x^2}\phi (z,x)
+z\frac{\pt}{\pt x}\phi (z,x)}
\right]\\
\qquad\qquad\qquad\qquad\qquad+\sigma (z)e^{-2zx}\int_\bbR T^{(x)}(z,\rd z^\prime )[\phi (x,z^\prime)
-\phi (x,z)],\ena \quad x,z\in\bbR .\eeq     
$\Big($A simple example 
is where $z,z^\prime =\pm 1$ with jump intensities $T^{(x)}(1,-1)=T^{(x)}(-1,1)=q(x)$; here the IM 
$\ups =b(\delta_{1}+
\delta_{-1})$, the sum of Dirac deltas supported at $z=\pm 1$, with equal coefficients. Another example is 
where $T^{(x)}(z,\rd z^\prime )=T^{(x)}(z+u,\rd (z^\prime +u))$ \ $\forall$ $u,z,z^\prime\in\bbR$ 
(shift-invariant jump measures). Here $\ups (\rd z)=\rd z$ is Lebesgue.$\Big)$
\vskip .5 truecm

{\bf Theorem 5.4.} (I) {\sl Operator $\BR$ and measure $\bka$ satisfy the WIE on $\bbR\times\bbR$.} (II) 
{\sl Under the SCC $\Xi:=\int_{\bbR\times\bbR}{\diy\frac{e^{2zx}}{\sigma (z)}}\rd x\ups (\rd z)<\infty$, we 
obtain a unique PD satisfying the WIE.}
\vskip .5 truecm

{\bf Remarks.\ 5.3.} As above, one can refer to Theorem 3.1, P. 377 in \cite{EK}, and provide 
sufficient conditions upon $T^{(x)} (z,\rd z^\prime )$ under which the combined MP will exist, with 
a Feller transition semigroup on $\bbR\times\bbR$.

{\bf 5.4.} The SCC $\Xi <\infty$ requires a rapid growth of $\sigma (z)$ at $z\sim\pm\infty$.  
\vskip .5 truecm

{\bf 5.C. Models of a two-component diffusion.} The next class of models is where both basic and SE 
MPs are diffusions. The simplest model is where the basic MP ${\wt X}^{(z)}$ is a $d$-dimensional 
time-scaled Wiener process (WP), with 
$\bbX=\bbR^d$, the SE space $\cZ=\bbR$, and the SE process ${\wt Z}(t)$  is a standard WP on a line.
Here 
\beq\label{eq:bscWnrL}\beac \tL^{(z)}g(x)=\diy\frac{z^2}{2}\sum_{1\leq i\leq d}
\frac{\pt^2}{\pt\rx_i^2} g(x),\\ 
\tA f(z)=bf^\prime (z)+\diy\frac{1}{2}f^{\prime\prime}(z),\ena
\quad x=(\rx_1,\ldots ,\rx_d)\in\bbR^d,\;z\in\bbR.\eeq
Moreover, the reference measure $\rd\gam (x)$ is $d$-dimensional Lebesgue, $\rd x$, and  
the IM $\rd \ups (z)$ is absolutely contiunous wrt $\rd z$, one-dimensional Lebesgue: 
\beq\label{eq:bkaW}\beac m(z,x)=1\;\hbox{ that is, $\rd\nu^{(z)}(x)=\rd x$,}\\
\rd\ups (z)=e^{2bz}\rd z .\ena\eeq 
Both $\nu^{(z)}$ and $\ups$ are genuine IMs for their respective MPs ${\wt X}^{(z)}$ and ${\wt Z}(t)$. 

Next, fix coefficient functions $z\in\bbR_+\mapsto\alpha (z)>0$ and $z\in\bbR_+\mapsto\sigma (z)>0$
such that 
\beq\label{eq:uslo}\hbox{$\alpha$ and $\sigma$ are Hoelder functions, and $c<z^2\alpha (z),\sigma (z)<C(1+z^2)$}\eeq
where $c,C\in (0,\infty )$ are constants. The combined
generator $\BR$ has the form
\beq\label{eq:BMdffsn}\BR\phi (z,x)=\alpha (z)\frac{z^2}{2}\sum_{1\leq i\leq d}
\frac{\pt^2}{\pt\rx_i^2}\phi (z,x)+\sigma (z)\left[b\frac{\pt}{\pt z}\phi (z,x)+\frac{1}{2}\frac{\pt^2}{
\pt z^2}\phi (z,x)\right],\;x\in\bbR^d,\,z\in\bbR,\eeq
and acts on $\rC^2$-functions $(z,x)\mapsto\phi (z,x)$. More precisely, $\BR$ is a closed operator 
whose domain $\rD (\BR)$ is the closure, in $\rC_\rb (\bbR_+\times\bbR^d)$,  of the set of  functions 
$\phi$ such that $\BR\phi\in\rC_\rb (\obR_+\times\bbR^d)$.
Owing to a general result, such as Theorem 1.5 on P. 369 in  \cite{EK},  Section 8.1; see also Ref. \cite{KS}, 
Sections 5.2--5.4, $\BR$ gives rise to a unique Feller semi-group of operators described by 
transition densities (relative to $\rd z\times\rd x$). Consequently, 
there exists a unique two-component diffusion process $(Z(t),X(t))$ in $\bbR_+\times\bbR^d$
generated by $\BR$. Then, in accordance with Proposition 9.2 in \cite{EK}, measure $\bka$ with
\beq\label{eq:bkadffsn}\bka (\rd z\times\rd x)=\frac{e^{2bz}}{\sigma (z)}\rd z\rd x\eeq
yields a genuine IM for $(Z(t),X(t))$. As the basic IM is Lebesgue,  $\Xi :=\int_{\bbR_+\times\bbR}\diy
\frac{e^{2bz}}{\sigma (z)}\rd z\rd x=\infty$, and there is no EPD.

The two-component diffusion can be described as a solution to
a system of stochastic integro-differential equations  (SIDEs)  
\beq\label{eq:comdif1}\beac 
\rd X(t)= \sqrt{ \alpha(Z(t))}Z(t)\rd W_{1}(t),\; X(0)=x^0\in\bbR^d,\\
\rd Z(t)=  b\sigma (Z(t))dt+ \sqrt{\sigma(Z(t))}\rd W_{2}(t) +\rd L_{Z}(t),\;Z(0)=z^0>0,\\
L_{Z}(t)=\int_{0}^{t}{\mathbf 1}(Z(s)=0)\rd L_{Z}(s).
\ena \eeq
Here $W_1(t)$ is the WP in $\bbR^d$ and $W_2(t)$ the WP in $\bbR$, and 
$W_1(t)$ and $W_2(t)$ are independent. Process $L_{Z}(t)$ is a local  time spend by 
$Z(s)$ at $0$ by time $t$. ($L(t)$ increases on a set of time points of measure $0$.)
\vskip .3 truecm

{\bf Remark 5.4.} SDEs \eqref {eq:comdif1} involve the solution to the so-called Skorohod
problem. Here we consider a WP  $W(t)$ on $\bbR$  and construct a BM $Z(t)$ on $(0,\infty )$
with reflection at $0$ by taking a pair of processes
 $(Z(t),L(t))$ such that (i) $Z(t)=z^0+W(t)+L(t)$, $t\geq 0$, $z^{0}>0$, (ii) $Z(t)\ge 0$, $t\geq 0$,  
 and (iii) process $L(t)$ has $L(0)=0$ (starts at $0$), is continuous, monotone nondecreasing and obeys
 ${\rm{supp}}\,(\rd L)\subset\{t\ge 0:Z(t)=0\}$  (that is, $L(t)$ increases only at times when $Z(t)$ is equal to 
 zero).  Moreover, (iv) $L(t)=[-z^0-\inf_{s\le t}W(s)]^+$. In addition,  (v) the distributions of processes 
 $Z(t)$  and  $|W(t)|$  coincide.  
 
 In a similar way, for ${\cal Z}=(0,1)$ the solution of the Skorohod problem is a triple of processes
 $(Z(t),L(t),U(t))$ such that (i) $Z(t)=z^0+W(t)+L(t)-U(t)$, $t\geq 0$, $0<z^{0}<1$, (ii) $0\leq Z(t)\leq 1$, 
 $t\geq 0$,  (iii) process $L(t)$ 
has $L(0)=0$, is continuous, monotone nondecreasing and obeys
 ${\rm{supp}}\,(\rd L)\subset\{t\ge 0:Z(t)=0\}$, and (iv) a similar property holds for $U(t)$ (repacing 
 endpoint $0$ with $1$).  

An analogous construction works for $\cZ =(z^1,z^2)$. 

\section{Two-component diffusions based on an Ornstein--Uhlenbeck process}\label{sec:6}

In this section we consider several examples of combined two-dimensional diffusion
where the basic MP is an 
Ornstein--Uhlenbeck process (OUP). Such examples may be of interest in financial calculus.

In all models 
that are discussed here, (a) we take $\alpha (z)=\sigma (z)\equiv 1$, (b) the 
SE generator $\tA^{(x)}=\tA$ is independent of $x$, (c) 
space $\bbX$ is a line $\bbR$ or a bounded interval $\bbI=(x^1,x^2)\subset\bbR$, with Lebesgue 
measure $\rd\gam (x)=\rd x$, and (d) space $\cZ$ coincides with a half-line $\bbR_+=(0,\infty )$ 
or the whole line $\bbR$ or a bounded interval $\bbJ=(z^1,z^2)\subset\bbR$, and (e) $\ups$ is 
an absolutely continuous measure $\ups$: $\rd\ups (z)=w(z)\rd z$. (We will loosely refer to $\cZ$ as an
interval.)  The basic MP ${\wt X}^{(z)}(t)$ on $\bbX$ has generator $\tL^{(z)}$ where, for a $\rC^2$-function
$x\in\bbX\mapsto g(x)$, 
\beq\label{eq:L501}
\tL^{(z)}g(x)=-xg^\prime (x)+\frac{z^2}{2}g^{\prime\prime}(x),\;\hbox{ with }\;
m(z,x)=e^{-x^2/z^2},\;\;x\in\bbX,\;z\in\cZ.
\eeq
 More precisely, let $\rC_\rb=\rC_\rb (\bbX)$ denote the space of bounded continuous functions 
 $x\in\bbX\mapsto g(x)$
with the sup-norm. The (closed) 
operator $\tL^{(z)}$ acts on its domain $\rD (\tL^{(z)})\subset\rC_\rb(\bbX)$ which is  the closure of the set of 
 $\rC^2$ functions $g$ such $\tL^{(z)}g\in\rC_\rb (\bbX)$ and $g$ satisfies a chosen boundary condition 
 at the endpoints in case $\bbX=\bbI$. 
The coefficient $z^2$ is interpreted as a stochastic volatility. In all models under consideration, there 
exists a Feller semi-group generated by $\tL^{(z)}$ which is determined by transition probability 
densities. Viz., for $\bbX=\bbR$,  the densities are
\beq\label{eq:trdnsOU}\diy p_t(x,x^\prime)=\frac{1}{{\sqrt\pi}|z|}\exp\,
\left[-\frac{(x^\prime-xe^{-t})^2}{z^2(1-e^{-2t})}\right],\;\;x,x^\prime\in\bbR.\eeq 
In this case, after normalization, the IM density $m(z,x)=e^{-x^2/z^2}$, $z\neq 0$, gives rise to a Gaussian 
probability distribution N$(0,z^2/2)$.  

Examples of the SE diffusion in this section  will be a Brownian motion (BM) on $\bbR_+$ with a 
reflection at $0$ (sub-section 6.B), an OUP
on $\bbR$ (sub-section 6.C) and an affine Cox--Ingersoll--Ross diffusion on $\bbR_+$  
(sub-section 6.D), as well as versions of these models on bounded intervals $\bbJ$. 

The existence/uniqueness of combined Feller MPs on bounded rectangles $\bbI\times\bbJ$ will follow 
from regularity (Hoelderness and non-degeneracy) of coefficients/boundary conditions.  See, Ref. 
\cite{EK},  Section 8.1, Theorem 1.5 on P. 369, or Ref. \cite{KS}, Sections 5.2--5.4. Consequently, for 
such examples covered by these general results we 
will be able to use the term IM unreservedly. Otherwise we only work with generator $\BR$ (see
Eqn \eqref{eq:combGen}) and will have to employ the WIEs (weak-invariance equations). 
Cf. Eqn \eqref{eq:diffWIE} below.
\vskip .3 truecm

{\bf 6.A.} Let us first re-formulate our general construction in the situation where basic and 
environment  MPs
are diffusions. As was said, we consider the case where the dimension of each component equals $1$;
however, the general scheme can be also developed in a multi-dimensional setting. The attention is 
on a two-component diffusion MP $(Z (t),X(t))$ in a Cartesian product 
$\cZ\times\bbX$  where $\bbX$ and $\cZ$ are of the types described above.   
Here, the Lebesgue measures $\rd z$ and
$\rd x$ will play special roles.  

The coefficient functions
\beq\label{eq:cffcntHld} x\in\bbR\mapsto a(x)\in\bbR,\; z\in\cZ\mapsto c(z)\in\bbR,\; z\in\cZ\mapsto 
C(z)>0, \;
(z,x)\in\cZ\times\bbX\mapsto m(z,x)>0\eeq 
are supposed to be Hoelder,  and  $W_1(t)$ and $W_2(t)$ are 
independent standard Wiener processes (WPs) in $\bbR$. We also specify that,
for a diffusions with accessible boundaries we use the Neumann condition. 

Accordingly,  the generator $\BR$ is a closed operator in $\rC_\rb (\cZ\times\bbX )$;
on $\rC^2$-functions $(z,x)\in\cZ\times\bbX\mapsto\phi (z,x)$ its action is given by
\beq\label{eq:combGen}\beacl\BR\phi(z,x)&\diy=a(x)\,\frac{\pt \phi}{
\partial x}(z,x)+
\frac{z^2}{2}\frac{\pt ^{2}\phi}{\pt x^{2}}(z,x)\\
\;&\diy \qquad+\frac{1}{m(z,x)}\left[c(z)\,\frac{\pt \phi}{\pt z}(z,x)
+\frac{C(z)}{2}\frac{\pt^{2} \phi}{\pt z^{2}}(z,x)\right].\ena \eeq
More precisely, the domain $\rD(\BR)$ is the closure of $\rC^2$-functions $\phi$ for which  
$\BR\phi\in\rC_\rb(\ocZ\times\obX)$, and the Neumann boundary conditions are fulfilled 
on the accessible parts of the boundary.
Here $\ocZ$ and $\obX$ stand for the closure of $\cZ$ and $\bbX$
and $\rC_\rb(\ocZ\times\obX)$
denotes the space of bounded continuous functions on $\ocZ\times\obX$. 

Viz., for $\bbX=\bbR$ and $\cZ =\bbR_+=(0,\infty )$,  
assuming that $0$ is an instant reflection point for ${\wt Z}(t)$,
generator $\BR$ introduced in \eqref{eq:combGen} is defined on a domain $\rD (\BR)$ 
which is the closure of $\rD^0(\BR)$ where
\beq\label{eq:gengen}\beal\rD^0(\BR)=\Big\{\phi\in\rC^2(\bbR_+\times\bbR):\;\BR\phi\in
\rC_\rb (\obR_+\times\bbR),\\
\qquad\qquad\quad\diy
\;\BR\phi (0+,x)=\lim\limits_{z\to 0+}\BR\phi (z,x),\;\frac{\pt\phi (z,x)}{\pt z}\Big|_{z=0+}=0\;\;
\forall\;x\in\bbR\Big\}.\ena\eeq 

As was mentioned,
in some examples the existence and uniqueness of a Feller diffusion $(Z(t),X(t))$ follows from general 
existence/uniqueness theorems; when such known results are not applicable, the notation
$(Z(t),X (t))$ and the term a combined process have only an inspirational meaning, and the term 
IM is a euphemism for a solution to the WIE in Eqn \eqref{eq:diffWIE}.
 
\vskip 3 truemm

To construct an IM $\bka $ for $(Z(t),X (t))$, we make two assumptions. First, we suppose that, $\forall$ 
given $z$, function $x\mapsto m(z,x)$ gives the density of an IM $\nu^{(z)}$ for a diffusion MP 
${\wt X}^{(z)}(t)$ in
$\bbR$ with some boundary conditions. 
The generator $\tL^{(z)}$ of process ${\wt X}^{(z)}$ is a closed operator in $\rC_\rb(\bbX)$ 
whose action 
on $\rC^2$-functions $x\mapsto g(x)$ is given by \beq\label{eq:bscdffGen} \tL^{(z)}g(x)=
a(x)g^\prime (x)+\frac{z^2}{2}g^{\prime\prime}(x),\;\;x\in\bbX.\eeq
More precisely, the domain $\rD (\tL^{(z)})$ is the closure of the set of $\rC^2$-functions $g$ for which
$\tL^{(z)}g\in\rC_\rb (\obX)$ and the derivative $g^\prime (x)$ vanishes at the
accessible points of the boundary $\pt\bbX$.

In other words, we assume that for a sufficient amount of functions $g$ (forming a core for $\tL^{(z)}$)
we have that 
\beq\label{eq:bscdffIM}\int_\bbD m(z,x)\,\tL^{(z)}g(x)\rd x=0,\;\hbox{ implying }\;
{\tL^{(z)}}^*m(z,x)=0,\;\;z\in\cZ,x\in\bbX.\eeq
Here ${\tL^{(z)}}^*$ is the (properly defined) adjoint operator acting on Radon--Nikodym densities 
(relative to $\rd x$). The reader familiar with the concepts of the scale and speed densities can think
that 
$$m(z,x)\propto\exp\,\left[\frac{2}{z^2}\int^x a(y)\rd y\right].$$  

Second, we assume that there is an SE diffusion MP ${\wt Z}(t)$ in interval $\cZ$ 
which obeys the SDE
\beq\label{eq:envdffSDE}\rd{\wt Z}(t)=c({\wt Z}(t))\rd t+ 
[C({\wt Z}(t))]^{1/2}\rd{\wt W}_2 (t),\;\;{\wt Z}(0)=z_0,\eeq
and has an IM $\ups$. (As above, 
the SDE \eqref{eq:envdffSDE} is subject to a modification at point 
$0$ when $\cZ=(0,\infty )$ and $0$ is an accessible boundary.) Here ${\wt W}_2(t)$ is a standard Wiener process in $\bbR$,
and $\rd\ups (z)$ is assumed to be absolutely continuous wrt $\rd z$: $\rd\ups (z)=w(z)\rd z$. The 
generator $\tA$ of process ${\wt Z}(t)$ is a closed operator acting on  $\rC^2$-functions $z\mapsto f(z)$ by  
\beq\label{eq:envdffGen}\tA f(z)=c(z)f^\prime (z)+\frac{C(z)}{2}f^{\prime\prime}(z),\;\;\;z\in\cZ ;\eeq
its domain $\rD(\tA)$ is the closure of the set of $\rC^2$-functions $f$ with $\tA f\in\rC_\rb(\ocZ)$ 
such that $f^\prime =0$ at accessible endpoints of $\cZ$. Thus, it is assumed that  for sufficiently many functions 
$f$ (forming a core for $\tA$) we have that 
\beq\label{eq:envdffIM}\int_\bbB \tA f(z)w(z)\rd z=0
\;\hbox{ implying }\;{\tA}^*w(z)=0,\;\;z\in\bbB.\eeq
Referring, as before, to   the scale and speed densities, one can assume that 
$$w(z)\propto\exp\,\left[\int^z\frac{2c(y)}{C(y)}\rd y\right].$$

We say that a measure $\bfta (\rd z\times\rd x)$ on $\cZ\times\bbX$  satisfies the WIE 
 with generator $\BR$ if 
\beq\label{eq:diffWIE}\int_{\cZ\times\bbR}\BR\phi (z,x)\bfta (\rd z\times\rd x)=0\eeq
for any function $(z,x)\in\cZ\times\bbR\mapsto\phi (z,x)$ satisfying conditions (i)--(ii) below. 
(i) $\forall$ $z$,  the section map $g_{\phi ,z}:\;x\to\phi (z,x)$ 
lies in a core of operator $\tL^{(z)}$. (ii) $\forall$ $x$,  the section map $f_{\phi ,x}:\;z\to\phi (z,x)$ 
lies in a core of operator $\tA$. (iii) The following integrals are finite:
$$\int_{\cZ\times\bbR}\left|a(x)\frac{\pt \phi}{\pt x}(z,x)+
\frac{z^{2}} 2\frac{\pt ^{2}\phi}{\pt x^{2}}(z,x) \right|\bfta (\rd z\times\rd x) $$
and
$$\int_{\cZ\times\bbR}\frac{1}{m(z,x)}\left|c(z)\frac{\pt \phi}{\pt z}(z,x)
+\frac{C(z)}2\frac{\pt^{2} \phi}{\pt z^{2}}(z,x)\right|\bfta (\rd z\times\rd x),$$
which allows us to use any order of integration in the summands emerging in the LHS of 
\eqref{eq:diffWIE}.

\vskip .5 truecm

{\bf Theorem 6.1.} {\sl Assume \eqref{eq:bscdffIM} and \eqref{eq:envdffIM}. Then the following
assertions hold true.} (I) {\sl The measure $\bka$ with the  Radon--Nikodym density 
\beq\label{eq:cmbIM}
\diy\frac{\bka (\rd z\times \rd x)}{\ups (\rd z)\times\rd x}= m(z,x)\eeq
satisfies the WIE with generator $\BR$ from Eqn \eqref{eq:combGen}.} (II) {\sl 
In case of an absolutely continuous measure $\ups$, with $\ups (\rd z)=w(z)\rd z$ we obtain 
$\diy\frac{\bka (\rd z\times \rd x)}{\rd z\times\rd x}=m(z,x)w(z)$. } (III) {\sl Assume that $\bbX$ and $\cZ$ are
bounded intervals $\bbI\subset\bbR$ and $\bbJ\subset (\eps ,\infty )$, respectively, where $\eps >0$. 
Suppose that coefficient functions $a(x)$, $c(z)$, $C(z)$, $m(z,x)$ 
are bounded and Hoelder, with \ $\inf÷\;÷C(z)>0$ and \ $\inf÷\;÷m(z,x)>0$. Then there is a unique Feller 
diffusuion process $(Z(t),X(t))$ with generator $\BR$ 
in rectangle $\bbJ\times\bbI$, and measure $\bka$ is an IM for $(Z(t),X(t))$.}
\vskip .5 truecm

{\it Proof.} Let us start with (I). We have to check that 
$\int_{\cZ\times\bbR}\big[\BR\phi (z,x)\big]\bka (\rd z\times\rd x)=0$, for each 
function $(z,x)\mapsto\phi (z,x)$ mentioned in the above definition of the WIE. We write that 
$$\beal\diy\int_{\cZ\times\bbR}\BR\phi (z,x)\bka (\rd z\times\rd x)= \int_{\cZ\times\bbR}\left\{\left[
a(x)\frac{\pt \phi}{\pt x}(z,x)+
\frac{z^{2}} 2\frac{\pt ^{2}\phi}{\pt x^{2}}(z,x)\right]\right.\\
\diy\left.\qquad\qquad\qquad+\frac{1}{m(z,x)}\left[c(z)\frac{\pt \phi}{\pt z}(z,x)
+\frac{C(z)}2\frac{\pt^{2} \phi}{\pt z^{2}}(z,x)\right]\right\}m(z,x)\ups (\rd z)\rd x\ena$$
which is equal to a sum of integrals $I_1+I_2$ where
$$\bear\diy I_1=\int\left[
a(x)\frac{\pt \phi}{\pt x}(z,x)+
\frac{z^{2}} 2\frac{\pt ^{2}\phi}{\pt x^{2}}(z,x)\right]m(z,x)\rd \ups (z)\,\rd x\qquad{}\\
\diy =\int_\cZ \left\{\int_\bbD m(z,x)
\left[\tL^{(z)}\phi (z,\,\cdot\,)\right](x)\rd x\right\}\rd\ups (z)\ena$$
and
$$ I_2=\int\left[c(z)\frac{\pt \phi}{\pt z}(z,x)
+\frac{C(z)}2\frac{\pt^{2} \phi}{\pt z^{2}}(z,x)\right]\ups (\rd z)\rd x
=\int_\bbD\left[\int_\cZ \big[\tA\phi (\,\cdot\,,x)\big] (z)\ups (\rd z)\right]\rd x .$$
In both expressions, the inner integrals vanish: in $I_1$ it occurs due to \eqref{eq:bscdffIM}, 
and in $I_2$ by virtue of \eqref{eq:envdffIM}.

Assertion (II) is straightforward, whereas (III) follows from general results (see, e.g., \cite{EK},  Section 
8.1, Theorem 1.5 on P. 369).
\vskip .5 truecm


We now pass to examples of interest. 
\vskip 3 truemm

{\bf 6.B.} \ We start with
an example where $\cZ=\bbR_+=(0,\infty )$. Here, the environment MP ${\wt Z}(t)$ is 
represented by a BM in $\bbR_+$ with drift $b\in\bbR$ and reflection at $0$.  Recall: given a 
standard WP $W(t)$,
there exists  a solution to the Skorohod problem. That is, $\forall$ $z_0>0$, $\exists$ a pair of processes
 $Z(t),L(t)$ such that (i) $Z(t)=z_0+W(t)+L(t)$, $t\geq 0$, (ii) $Z(t)\ge 0$, $t\geq 0$,  and (iii) process $L(t)$ 
has $L(0)=0$ (starts at $0$), is continuous, monotone nondecreasing and obeys
 ${\rm{supp}}\,(\rd L)\subset\{t\ge 0:Z(t)=0\}$  (that is, $L(t)$ increases only at times when $Z(t)$ is equal to 
 zero).  Moreover, (iv) $L(t)=[-z_0-\inf_{s\le t}W(s)]^+$. In addition,  (v) the distributions of processes 
 $Z(t)$  and  $|W(t)|$  coincide. Formally, we have:
\beq\label{eq:L421}\beac\diy\tL^{(z)}g(x)=-xg^\prime (x)+\frac{z^2g^{\prime\prime}(x)}{2},\;\;
m(z,x)=e^{-x^2/z^2},\\
\diy\tA f(z)=bf^\prime (z)+\frac{f^{\prime\prime}(z)}{2},\quad\rd\ups (z)=e^{2bz}\rd z,\qquad\ena \;\,x\in\bbR ,z>0.
\eeq
The domain of operator $\tL^{(z)}$ is the closure in $\rC_\rb(\bbR)$ of $\rC^2$-functions 
$x\in\bbR\mapsto g(x)$ such that $\tL^{(z)}g\in\rC_\rb(\bbR)$. The domain of
$\tA$ is the closure in $\rC_\rb (\bbR_+)$ of $\rC^2$-functions $z\in\bbR_+\mapsto f(z)$ such that 
$f^\prime (0+)=0$ and $\tA f\in\rC_\rb(\bbR_+)$. Both $\tL^{(z)}$ and $\tA$ generate unique Feller 
semi-groups. For
$\tL^{(z)}$, the Feller semi-group is determined by the transition densities $p_t(x,x^\prime )$ from 
Eqn \eqref{eq:trdnsOU}. For $\tA$, the Feller semi-group is determined by the transition densities
\beq\label{eq:trdnsrefl}\beacl
p^{(t)}(z,z^\prime)&\diy =\frac{2be^{2bz}}{e^{2bz}-1}+\frac{2}{\pi}e^{b(z^\prime -z)-b^2t/2}\\
\;&\qquad\times\int_0^\infty\diy\frac{e^{-s^2t/2}}{s^2+b^2}\Big[s\cos\,(sz)+b\sin\,(sz)\Big]\Big[s\cos\,(sz^\prime)
+b\sin\,(sz^\prime ) \Big],\;\;z,z^\prime >0.\ena\eeq
Cf. \cite{Lin}, Eqns (28)--(29) (more precisely, a displayed equation between (28) and (29).

In this example, for $x\in\bbR$ and $z>0$ we have: 
\beq\label{eq:BM421}\beacl\BR\phi (z,x)&\diy=-x\frac{\pt}{\pt x}\phi (z,x)
+\frac{z^2}{2}\frac{\pt^2}{\pt x^2}\phi (z,x)\\
\;&\diy\qquad+e^{x^2/z^2}\left[b\frac{\pt}{\pt z}\phi (z,x)+
\frac{1}{2}\frac{\pt^2}{\pt z^2}\phi (z,x)\right],\ena\eeq
and 
\beq\label{eq:bka421}\rd\bka (z,x)=e^{2bz- x^2/z^2}\rd x\,\rd z.\eeq
The domain $\rD(\BR)$ of operator $\BR$ in \eqref{eq:BM421} is the closure of the set $\rD^0(\BR)$ 
from \eqref{eq:gengen}. 
Due to presence of function $e^{x^2/z^2}=[m(z,x)]^{-1}$, the existing general results do not guarantee
that operator $\BR$ generates a unique Feller semi-group in $\rC_\rb (\bbR_+\times\bbR)$. 
Nevertheless, Theorem 8 (I) holds, and $\bka$ satisfies the WIE with $\BR$.  

However, if we take $\bbX$ to be an interval $\bbI =(x^1,x^2)$ and $\cZ=(z^1,z^2)$, with
$-\infty <x^1<x^2<\infty$ and  
$0<z^1<z^2<\infty$, and  assume that both ${\wt X}^{(z)}$ and ${\wt Z}(t)$ are reflected at the endpoints 
$x^i$ and $z^i$ then the situation changes. We still define $\BR$ and $\bka$ via Eqns \eqref{eq:BM421} 
and \eqref{eq:bka421}, for $x\in\bbI$, $z\in\bbJ$, but the domain $\rD(\BR)$ involves the Neumann
boundary conditions $\diy\frac{\pt}{\pt x}\phi (z,x^i)=0$ and $\diy\frac{\pt}{\pt z}\phi (z^i,x)=0$. Under 
these modifications, operator 
$\BR$ generates a unique Feller semi-group in $\rC_\rb (\bbJ\times\bbI)$. Furthermore, by 
Proposition 9.2 in \cite{EK}, P. 239, $\bka$ is a genuine IM for this sem-group. 

The corresponding MP is a pair $(X(t),Z(t))$ plus an auxiliary quadruple $(L_X(t),U_X(t),L_Z(t),U_Z(t))$ 
solving a system of SIDEs 
$$\beac \rd X(t)=-X(t)\rd t +Z(t)\rd W_{1}(t)+\rd L_{X}(t)-\rd U_{X}(t),\;X(0)=x^0\in (x_{1},x_{2})\\
\rd Z(t)=-be^{\frac{X^{2}(t)}{Z^{2}(t)}}\rd t+e^{\frac{X^{2}(t)}{2Z^{2}(t)}}\rd 
W_{2}(t)+\rd L_{Z}(t)-\rd U_{Z}(t),\;\; Z(0)=z^0\in(z_{1}, z_{2}),\\
\int_{0}^{t}{\mathbf 1}({X(s)=x_{1}})\rd L_{X}(s)=L_{X}(t),\quad \int_{0}^{t}{\mathbf 1}({X(s)=x_{2}})\rd U_{X}(s)=U_{X}(t),\\
\int_{0}^{t}{\mathbf 1}({Z(s)=z_{1}})\rd L_{Z}(s)=L_{Z}(t),\quad \int_{0}^{t}{\mathbf 1}({Z(s)=z_{2}})\rd U_{Z}(s)
=U_{Z}(t),\ena$$
where $L_{X}(t),U_{X}(t)$ are the local times spend by the process $X(t)$ at points $x_{1}$   and $x_{2}$ 
and $L_{Z}(t),U_{Z}(t)$ are local times spend by the process $Z(t)$ at points $z_{1}$   and $z_{2}$. 

Here and below,  $x_0$, $z_0$ are initial values and  $W_1(t)$, $W_2(t)$ are 
independent WPs in $\bbR$. 
\vskip 3 truemm

{\bf 6.C.}   In this example we take $\cZ=\bbR$ and deal with an OUP ${\wt X}^{(z)}(t)$ where 
the diffusion coefficient $z$ follows its own OUP. Formally speaking, we set 
\beq\label{eq:L423}\beac\diy\tL^{(z)}g(x)=-xg^\prime (x)+\frac{z^2}{2}g^{\prime\prime}(x),\;\;
m(z,x)=e^{-x^2/z^2}, \\
\diy\tA f(z)=-zf^\prime (z)+\frac{1}{2}f^{\prime\prime}(z),\;\;\rd\ups (z)=
e^{-z^2}\rd z,\quad\ena\;\;\;\;z,x\in\bbR .\eeq
Here, the Feller semi-group generated by $\tA$ is determined by the transition densities similar to 
\eqref{eq:trdnsOU}:
\beq\label{eq:trdnsOU1}\diy p^{(t)}(z,z^\prime )=\frac{1}{{\sqrt\pi}}\exp\,\left[-\frac{(z^\prime
-ze^{-t})^2}{(1-e^{-2t})}\right],\;\;z,z^\prime\in\bbR.\eeq 

Then, for $z,x\in\bbR$,  
\beq\label{eq:BM423}\begin{array}{cl}\BR\phi (z,x)&\diy =-x\frac{\pt}{\pt x}\phi (z,x)
+\frac{z^2}{2}\frac{\pt^2}{\pt x^2}\phi (z,x)\\
\;&\diy\qquad\qquad +e^{x^2/z^2}\left[-z
\frac{\pt}{\pt z}\phi (z,x)+\frac{1}{2}\frac{\pt^2}{\pt z^2}\phi (z,x)\right] ,
\ena \eeq
and 
\beq\label{eq:bka423}\rd\bka (z,x)=e^{-z^2- x^2/z^2}\rd x\,\rd z.\eeq
The domain $\rD(\BR )$ is the closure, in $\rC_\rb(\bbR\times\bbR)$, of $\rC^2$-functions 
$(z,x)\in\bbR\times\bbR\mapsto \phi (z,x)$ for which $\BR\phi\in\rC_\rb(\bbR\times\bbR)$.
(The boundary condition in \eqref{eq:BM423} is omitted.) 
Again, the available general results do not allow us to conclude that there exists a unique 
Feller semi-group generated by $\BR$, but Theorem 8 (I) holds true.   

However, as before, we can take $\bbX =(x^1,x^2):=\bbI$ and $\cZ =(z^1,z^2):=\bbJ$ where 
$-\infty <x^1<x^2<\infty$, $-\infty <z^1<z^2<\infty$ and $0\not\in\bbJ$, and use Eqns \eqref{eq:BM423}
and \eqref{eq:bka423} for $z\in\bbJ$, $x\in\bbI$, adding the Neumann boundary conditions
$\diy\frac{\pt}{\pt x}\phi (z.x^i)=0$ and $\diy\frac{\pt}{\pt z}\phi (z^i.x)=0$ in the definition of 
$\rD(\BR)$. Then $\BR$ generates a unique Feller semi-group in $\rC (\bbJ\times\bbI)$, and $\bka$ is 
a genuine IM.  

The corresponding diffusion MP is a pair $(X(t),Z(t))$ from the solution to SIDEs
$$\beac \rd X(t)=-X(t)\rd t +Z(t)\rd W_{1}(t)+\rd L_{X}(t)-\rd U_{X}(t),\;X(0)=x^0\in (x_{1},x_{2}),\\
\rd Z(t)=-Z(t)e^{\frac{X^{2}(t)}{Z^{2}(t)}}\rd t+e^{\frac{X^{2}(t)}{2Z^{2}(t)}}\rd W_{2}(t)+\rd L_{Z}(t)- 
\rd U_{Z}(t),\; Z(0)=z^0\in(z_{1}, z_{2}),\\
\int_{0}^{t}{\mathbf 1}({X(s)=x_{1}})\rd L_{X}(s)=L_{X}(t),\quad \int_{0}^{t}{\mathbf 1}({X(s)=x_{2}})\rd U_{X}(s)=U_{X}(t),\\
\int_{0}^{t}{\mathbf 1}({Z(s)=z_{1}})\rd L_{Z}(s)=L_{Z}(t),\quad \int_{0}^{t}{\mathbf 1}({Z(s)=z_{2}})\rd U_{Z}(s)=U_{Z}(t).\ena$$
\vskip .5 truecm

{\bf 6.D.} Finally, consider the case where the diffusion coefficient $z$ follows a 
Cox--Ingersoll--Ross process on $\bbR_+$: see \cite{CIR}. This example can be considered as 
a modification of the 
Heston model of a stochastic volatility; cf. \cite{He}. Here we start with 
$\cZ=\bbR_+$ and generator $\tA$ of the SE diffusion ${\wt Z}(t)$ of the form
\beq\label{eq:CIR}\tA f(z)=a(b-z)f^\prime (z)+\frac{z}{2}f^{\prime\prime}(z),\;\;
\rd\ups (z)=z^{2ab -1}e^{-2az}\rd z,\;z>0,\eeq
where $a\geq 0,b>0$ are given parameters. The domain $\rD(\tA )$ consists of $\rC^2$-functions 
$z\in\bbR_+\mapsto f(z)$ such that $\tA f\in\rC_\rb (\bbR_+)$ and -- when $a> 1/2$ -- the right
derivative $f^\prime (0+)=0$. The Feller semigroup 
generated by $\tA$ has the transition density 
\beq\label{eq:trdnsCIR}p^{(t)}(z,z^\prime )=\diy c\Big[\exp\,\big(-cze^{-at}-cz^\prime\big)\Big]
\left(\frac{z^\prime e^{at}}{z}\right)^{q/2}I_q\left(2c\sqrt{zz^\prime e^{-at}}\right),\;\;z,z^\prime >0,\eeq
where $c=\diy \frac{2a}{1-e^{-at}}$, $q=2ab-1$ and $I_q$ is the Bessel function of order $q$. (For $a\leq 1/2$, 
process ${\wt Z}(t)$ does not hit $0$, whereas for $a>1/2$ it hits $0$ at infinitely many times that are 
indefinitely large.) 

The basic process is, as before, an OUP, with generator $\tL^{(z)}$, as in \eqref{eq:L501}, 
\eqref{eq:L421} and \eqref{eq:L423}.
For $x\in\bbR$, $z>0$, the combined generator is
\beq\label{eq:HestMoM}\begin{array}{cl}\BR\phi (z,x)&\diy =-x\frac{\pt}{\pt x}\phi (z,x)
+\frac{z^2}{2}\frac{\pt^2}{\pt x^2}\phi (z,x)\\
\;&\diy\qquad +e^{x^2/z^2}\left[a(b-z)\frac{\pt}{\pt z}\phi (z,x)
+\frac{z}{2}\frac{\pt^2}{\pt z^2}\phi (z,x)\right],\ena \eeq
and we set
\beq\label{eq:HestMoka}\rd\bka (z,x)=z^{2ab -1}e^{-2az- x^2/z^2}
\rd x\,\rd z.\eeq
The domain $\rD(\BR )$ is the closure of $\rD^0(\BR)$ where
\beq\label{eq:dffnDH}\bear
\rD^0(\BR)=\Big\{\phi\in\rC^2(\bbR_+\times\bbR):\;\; 
\BR\phi\in\rC_\rb(\obR_+\times\bbR)
\hbox{ and -- for $a>1/2$ -- also}\qquad{}\\
\diy\frac{\pt}{\pt z}\phi (z,x)
\Big|_{z=0+}=0\;\hbox{ and }\;\lim\limits_{z\to 0+}\BR\phi (z,x)=\BR\phi (0+,x)\Big\};\ena\eeq 
cf. \eqref{eq:gengen}).  Again we can$^\prime$t derive that there is a unique Feller 
semi-group generated by $\BR$ but Theorem 8 (I) holds true.

Taking $\bbX =(x^1,x^2)$ and $\cZ =(z^1,z^2)$ with $-\infty <x^1<x^2<\infty$ and $0<z^1<z^2<\infty$ 
and introducing Neumann boundary conditions 
leads, as above, to a unique Feller semi-group for which $\bka$ is a genuine IM.    
The coresponding diffusion MP is a pair $(X(t),Z(t))$ solving, together with $(L_X(t),U_X(t),L_Z(t),U_Z(t))$,
the system 
$$\beac \rd X(t)=X(t)\rd t +Z(t)\rd W_{1}(t)+\rd L_{X}(t)-\rd U_{X}(t),\;X(0)=x^0\in (x_{1},x_{2})\\
\rd Z(t)=a(b-Z(t))e^{\frac{X^{2}(t)}{Z^{2}(t)}}\rd t+ \sqrt{Z(t)}e^{\frac{X^{2}(t)}{2Z^{2}(t)}}\rd W_{2}(t)+\rd 
L_{Z}(t)-\rd U_{Z}(t),\; Z(0)=z^0\in(z_{1}, z_{2}),\\
\int_{0}^{t}{\mathbf 1}({X(s)=x_{1}})\rd L_{X}(s)=L_{X}(t),\quad \int_{0}^{t}{\mathbf 1}({X(s)=x_{2}})\rd
U_{X}(s)=U_{X}(t),\\
\int_{0}^{t}{\mathbf 1}({Z(s)=z_{1}})\rd L_{Z}(s)=L_{Z}(t),\quad \int_{0}^{t}{\mathbf 1}({Z(s)=z_{2}})\rd U_{Z}(s)=U_{Z}(t).\ena$$
\vskip .3 truecm

A similar construction can be done if we consider the SE process ${\wt Z}(t)$ on an interval $(0,r)$
instead of $\bbR_+$; cf. Section 6.1 in Ref. \cite{Lin}.
\vskip .5 truecm

{\bf Remarks.\  6.1.} One of intriguing problems emerging from the above constuction (in its general form 
as well as in specific examples) is to solve (and in fact, to pose in a correct manner) an inverse problem.
(We can speak of a hidden Markov model.)
Informally, we ask: is it possible to represent  a given a random process $X(t)$ as a projection 
$(Z(t),X(t))\mapsto X(t)$ 
of a combined MP $(Z(t),X(t))$ obtained by means of the construction, where component $Z(t)$ describes
a dynamics of the SE? E.g., in the
context of sub-section 4.A: can one represent a birth-death process on $\bbZ_+$ 
as a result of a projection $(z,n)\in\cZ\times\bbZ_+\mapsto n$, from an MP on $\cZ\times\bbZ_+$
constructed in the above fashion? A similar question in the context of sub-sections 4.B is:
given a process with continuous paths on (a domain $\bbD$ in) $\bbR^d$, is it possible to represent 
it as a projection of an MP  in a higher dimension, obtained by means of the above construction? 
In terms of Section 5, when it is possible to represent a given process with continuous paths on 
(a domain $\bbD$ in) $\bbR^d$ as a projection of a  diffusion in higher dimension?
\vskip .3 truecm

{\bf 6.2.} A feature of he construction presented in the paper is that the factor $m(z,x)$ (and also 
$\sigma (z)$ in previous parts of the paper) enter both the
generator $\BR$ and IM $\bka$ in a mutually inverse fashion. It creates an interesting pattern: if 
$m(z,x)$ is small then we have growing coefficients in $\BR$ causing difficulties with proving the 
existence of the combined MP $(Z(t),X(t))$. But then we get an IM $\bka$ which may be made an 
EPD (at least, at a formal level of the corresponding WIE). This suggests that for the models under 
consideration some new kind of existence theorems are possible, under weaker assumptions 
than in \eqref{eq:uslo}. This may present an interesting direction for future development. 

\vskip 1 truecm

{\bf Acknowledgements}.  YB works under Grant RFBR 15-01-01453 and expresses her 
gratitude to the RFBR. YMS thanks Math Dept, PSU, and IHES, for support and hospitality. 
YMS thanks C. Burdzy,  G. Pang, E. Pechersky and A. Yambartsev for stimulating discussions.

\vskip .2 truecm

\noindent Y Suhov: \quad yms@statslab.cam.ac.uk; ims14@psu.edu

\end{document}